\documentclass{amsart}
\input cyrfont.def
\input cyracc.def
\newfont{\msbm}{msbm10 scaled 1200}
\newtheorem{theorem}{Theorem}[section]
\newtheorem{lemma}[theorem]{Lemma}
\newtheorem{proposition}[theorem]{Proposition}
\theoremstyle{definition}
\newtheorem{definition}[theorem]{Definition}

\newtheorem{Pb}[theorem]{Problem}
\theoremstyle{remark}
\newtheorem{remark}[theorem]{Remark}
\numberwithin{equation}{section}
\DeclareMathOperator{\im}{Im}
\DeclareMathOperator{\nul}{Nul}

\DeclareMathOperator{\rank}{rank}
\DeclareMathOperator{\trace}{trace}
\DeclareMathOperator{\diag}{diag}
\begin{document}

\title[Isoprincipal Deformations of
Rational Matrix Functions II]{Rational Solutions of the
Schlesinger \\
System and Isoprincipal Deformations of \\
Rational Matrix Functions II}

\author{Victor Katsnelson}
\address{Department of Mathematics, Weizmann Institute of Science,
Rehovot 76100, Israel}
\email{victor.katsnelson@weizmann.ac.il}
\thanks{Victor Katsnelson  thanks  Ruth and Silvia Shoham
for endowing the chair that supports  his respective research.
Victor Katsnelson was also supported by the Minerva foundation.}
\author{Dan Volok}
\address{Department of Mathematics,
Ben--Gurion University of the Negev, Beer-Sheva 84105, Israel}
\email{volok@math.bgu.ac.il}

\keywords{Isoprincipal, isosemiresidual, joint system representation, Fuchsian
system, Schlesinger system} \subjclass{Primary: 47A56; Secondary: 34M55.}
\begin{abstract} In this second article of the series
we study holomorphic families
 of generic rational matrix functions parameterized
by the pole and zero loci. In particular, the isoprincipal deformations of generic
rational matrix functions are proved  to be isosemiresidual. The corresponding
 rational solutions of the Schlesinger
system are constructed and  the explicit expression for the related tau function is
given. The main tool is the theory of joint system representations for rational
matrix functions with prescribed pole and zero structures.
\end{abstract}
\maketitle

\section*{NOTATION}
\begin{itemize}
\item
 \(\mathbb C\)  stands for the complex plane.

\item
 \(\mathbb C_*\)  stands for the punctured complex plane:
\[\mathbb C_*=\mathbb C\setminus\{0\}.\]

\item
 \(\overline{\mathbb C}\)  stands for the extended
complex plane (\(=\) the Riemann sphere):
\[\overline{\mathbb C}=\mathbb{C}\cup\infty.\]

\item
\(z\) stands for the complex variable.

\item
 \({\mathbb C}^n\)  stands for the
\(n\)-dimensional complex
space.

\item
 In the coordinate notation, a point
\(\boldsymbol{t}\in{\mathbb
C}^n\) will be written as \(\boldsymbol{t}=(t_1,\ldots,t_n).\)

\item
\(\mathbb{C}^n_*\) is the set of
points \(\boldsymbol{t}\in{\mathbb C}^n,\)
whose coordinates \(t_1,\dots,t_n\)  are pairwise different:
\[{\mathbb{C}}^n_{\ast}={\mathbb C}^n\setminus\bigcup_{\substack{1\leq
i,j\leq n\\ i\not=j}} \{\boldsymbol{t}:t_i=t_j\}.\]

\item \(\mathbb{C}^{m\times n}\) stands for
the set of all \(m\times n\) matrices with complex entries.

\item For
 \(A\in\mathbb{C}^{m\times n},\) \(A^*\in\mathbb C^{n\times m}\) is the
adjoint matrix, \(\im(A)\) is the image
 subspace of \(A\) in \(\mathbb C^m\) (= the
linear span of the  columns of \(A\)) and \(\nul(A)\) is the null subspace of \(A\)
in \(\mathbb C^n\).

\item
 \([\cdot,\cdot]\) denotes the commutator: for
 \(A,B\in\mathbb{C}^{m\times m},\ [A,B]=AB-BA\).

\item\(I\)  stands for the identity
matrix of an appropriate dimension.

\end{itemize}

\setcounter{section}{9}
\section{Simple singularity of a meromorphic matrix function}
\label{GSMMF}

(For Sections 1\,-\,9 see \cite{KaVo}, \cite{KaVo1-e}.)

\begin{definition}%
\label{DefSimPol}%
Let \(R(z)\) be a \(\mathbb{C}^{m\times m}\)-valued
 function,  holomorphic in a
punctured neighborhood of a point \(t\in\mathbb C\). The point \(t\) is said to be a
{\em simple pole} of the matrix function \(R\) if
\begin{equation*}
R(z)=\frac{R_t}{z-t}+ H(z),
\end{equation*}
where \(R_t\in\mathbb{C}^{m\times m}\) is a constant matrix and the  function \(H\)
is holomorphic at the point \(t\). The matrix \(R_t\) is said to be the {\em
residue}  of the  function \(R\) at the point \(t\). Furthermore, if
\(r=\rank(R_t)\) and \(f_t\in\mathbb{C}^{m\times r}\) and
\(g_t\in\mathbb{C}^{r\times m}\) are matrices providing the factorization
\(R_t=f_tg_t,\) we shall say that \(f_t\) is the {\em left semiresidue} of \(R\) at
\(t\) and \(g_t\) is the {\em right semiresidue} of \(R\) at \(t\).
\end{definition}
\begin{remark}\label{gause}
The left and right semiresidues \(f_t,g_t\) are defined up to the transformation
\[f_t\rightarrow f_tc,\quad g_t\rightarrow c^{-1}g_t,\]
where  \(c\in\mathbb{C}^{r\times r}\) is an  invertible matrix.
\end{remark}
\begin{definition}%
\label{DefRegSingPoint} Let  \(R(z)\) be a \(\mathbb{C}^{m\times m}\)-valued
 function, holomorphic and invertible
 in a
punctured neighborhood of a point \(t\in\mathbb C\).
\begin{enumerate}
\item The point \(t\) is said to be {\em regular} for the function \(R\) if both the
function \(R\) and the inverse function \(R^{-1}\) are holomorphic functions in the
entire (non-punctured) neighborhood of the point \(t\), i.e. if \(R\) and \(R^{-1}\)
are holomorphic at the point \(t\).

 \item The point \(t\) is said to be {\em singular} for the function
\(R\) if at least one of the functions \(R\) and \(R^{-1}\) is
not holomorphic at
the point \(t\).
\end{enumerate}%
In particular, the point \(t\) is singular for the function \(R\) if \(R\) is
holomorphic at the point \(t\), but its value \(R(t)\) is a degenerate matrix. In
this case, the point \(t\) is said to be a {\em zero} of the function \(R\).
\end{definition}

\begin{definition}%
\label{DefGenSingPoint}%
Let  \(R(z)\) be a \(\mathbb{C}^{m\times m}\)-valued
 function, holomorphic and invertible
 in a
punctured neighborhood of a point \(t\in\mathbb C\), and let \(t\) be a singular
point of \(R\).
 The singular point  \(t\)  is said to be {\em simple} if one of
the following holds:

\begin{enumerate}

\item The point \(t\) is a  simple pole of the  function  \(R\)
 and a holomorphy point of the inverse function \(R^{-1}\).

 \item The point \(t\) is a  simple pole of the inverse function
\(R^{-1}\) and  a holomorphy point of the  function  \(R\) itself.

\end{enumerate}
\end{definition}%

\begin{remark}
Note that, according to Definiiton \ref{DefGenSingPoint}, if \(t\) is a simple
singular point of the function \(R\) then \(R\) is a single-valued meromorphic
function in the entire (non-punctured) neighborhood of \(t\).
\end{remark}

Our main goal is to study a matrix function  in a neighborhood of its simple
singular point from the point of view of linear differential systems. Thus we
consider the left logarithmic derivative of the function \(R\):

\begin{equation*}%
Q^l_R(z)\stackrel{\textup{\tiny def}}{=}R^{\prime}(z)R(z)^{-1}.
\end{equation*}%

\begin{remark}\label{leftder}
One can also consider the right logarithmic derivative of  \(R\):
\[
Q^r_R(z)=R(z)^{-1}R^{\prime}(z).
\]
But then \(-Q^r_R\) is the left logarithmic derivative of the
inverse function \(R^{-1}\):
\[Q^l_{R^{-1}}(z)=(R^{-1}(z))'R(z)=-R(z)^{-1}R'(z)R(z)^{-1}R(z)
=-Q^r_R(z).\] Thus in this work we shall deal mainly with left
logarithmic derivatives. Therefore, we shall use the notation
\(Q_R\) instead of \(Q^l_R\):
\[Q_R(z)\stackrel{\textup{\tiny def}}{=}R^{\prime}(z)R(z)^{-1},\]
and omit the word "left" when referring to the left logarithmic derivative.
\end{remark}


\begin{proposition}
\label{reslogder}
Let  \(R(z)\) be a \(\mathbb{C}^{m\times m}\)-valued
 function, holomorphic and invertible
 in a
punctured neighborhood of a point \(t\in\mathbb C\), and  let \(t\) be a simple
singular point of \(R\).
 Then the point \(t\) is a simple pole for  the  logarithmic
derivative\footnote{See Remark \ref{leftder}.} \(Q_R\) of \(R\). Moreover, for the
residue and the constant term of the Laurent expansion
\begin{equation}\label{lormr}
 Q_R(z)=\frac{Q_t}{z-t}+C+o(1)\text{ as }z\to t\end{equation}
the following relations hold.
\begin{enumerate}
\item If \(t\) is a pole of \(R\)  then
\begin{subequations}\label{necpo}
\begin{align}\label{necpo1}
 Q_t^2&=-Q_t,\\
\label{necpo3} Q_tCQ_t&=-CQ_t
\end{align}
\end{subequations}
and
\begin{equation}
 \label{necpo2} \im(Q_t)=\im(R_t),
 \end{equation}
where \(R_t\) is the residue of \(R\) at \(t\).

\item If \(t\) is a zero of \(R\)  then
\begin{subequations}\label{necze}
\begin{align}\label{necze1}
 Q_t^2&=Q_t, \\ \label{necze3}
Q_tCQ_t &=Q_tC\end{align}\end{subequations} and
\begin{equation}
 \nul(Q_t)=\nul(R_t),\label{necze2}
\end{equation}
where \(R_t\) is the residue of \(R^{-1}\) at \(t\).
\end{enumerate}
\end{proposition}
\begin{proof}
First, let us assume that \(t\) is a pole
of \(R\) and let
\begin{align}
R(z)&=\frac{R_t}{z-t}+A_0+A_1(z-t)
)+A_2(z-t)^2 +\ldots, \label{lorR}\\
R^{-1}(z)&=B_0+B_1 (z-t)+B_2(z-t)^2 + \ldots,
\label{lorR^{-1}}
\end{align}
be the Laurent expansions of the functions \(R\) and \(R^{-1}\) at
\(t\). Then
\begin{equation}
R^{\prime}(z)=-\frac{R_t}{(z-t)^2}+ A_1 + 2A_2(z-t)+ \ldots \label{deri}
\end{equation}
Multiplying the Laurent expansions term by term, we obtain from
 \eqref{lorR^{-1}} and \eqref{deri}
\begin{equation}
Q_R(z)= -\frac{R_tB_0}{(z-t)^2} -\frac{R_tB_1}{z-t} -R_tB_2+A_1B_0 +o(1).
\label{Rlder}
\end{equation}
Substituting the  expansions \eqref{lorR},
 \eqref{lorR^{-1}} into the identity
\[R^{-1}(z)R(z)=R(z)R^{-1}(z)=I,\]
we observe that
\begin{equation}\label{later1}
 R_tB_0=B_0R_t=0\text{ and } R_tB_1 + A_0B_0=I.
\end{equation}
Hence the first term of the expansion \eqref{Rlder} vanishes
and we obtain the
expansion \eqref{lormr} with
\begin{align}%
\label{Rres}%
Q_t&=-R_tB_1=A_0B_0-I,\\
\label{LEECT}%
C&=-R_tB_2+A_1B_0.
\end{align}%

Thus
\begin{equation*}%
(I+Q_t)Q_t= (A_0B_0) (-R_tB_1)= -A_0(B_0R_t)B_1=0,
\end{equation*}%
 i.e. \eqref{necpo1} holds. Furthermore,
\begin{multline*}%
(I+Q_t)CQ_t=
(A_0B_0)(-R_tB_2+A_1B_0)(-R_tB_1)=\\
=A_0(B_0R_t)B_2R_tB_1-A_0B_0A_1(B_0R_t)B_1=0,
\end{multline*}%
i.e. \eqref{necpo3} holds as well.
 Finally,
\begin{equation*}%
Q_tR_t=(A_0B_0-I)R_t= A_0(B_0R_t)-R_t=-R_t,
\end{equation*}%
which, together with \eqref{Rres}, implies \eqref{necpo2}.
 This completes  the proof
in the case when \(t\) is a pole of \(R\).
The case when \(t\) is a zero of \(R\)
can be treated  analogously.
\end{proof}
\begin{remark}\label{semi1}
Since for any \(p\times q\) matrix \(A\) the subspace \(\nul(A)\) is the orthogonal
complement of the subspace \(\im(A^*)\) in \(\mathbb{C}^q\), the relation
\eqref{necze2} can be rewritten as
\[\im(Q_t^*)=\im(R_t^*).\]
The latter relation, together with \eqref{necze1}, means that \(Q_t^*\) is a
(non-orthogonal, in general) projector onto the subspace
\(\im(R_t^*)\subset\mathbb{C}^m.\) Hence the right semiresidue of \(R^{-1}\) at its
pole \(t\) is also the  right semiresidue of \(Q_R\) at \(t\).

 Analogously, the
relations \eqref{necpo1} and \eqref{necpo2} mean that \(-Q_t\) is a (non-orthogonal,
in general) projector onto the subspace \(\im(R_t)\subset\mathbb{C}^m.\) Hence the
left semiresidue of \(R\) at its pole \(t\) is also the  left semiresidue of \(Q_R\)
at \(t\).
\end{remark}

 Proposition \ref{reslogder}  implies that
 a \(\mathbb{C}^{m\times m}\)-valued
function \(R(z)\) in a punctured neighborhood of its simple singular point \(t\) may
be viewed as a fundamental solution of a linear differential system
\begin{equation}\label{Fueq}
R^\prime(z)=Q(z)R(z),
\end{equation}
for which \(t\) is a {\em Fuchsian singularity} (see the first part of this work
\cite{KaVo} for details and references) and whose coefficients satisfy the relations
\eqref{necpo} or \eqref{necze}. The next proposition shows that \eqref{necpo} or
\eqref{necze} are the only requirements a differential system \eqref{Fueq} with a
Fuchsian singularity \(t\) has to satisfy in order for its fundamental solution in a
punctured neighborhood of \(t\) to be single-valued and have a simple singular point
at \(t\):
\begin{proposition}%
\label{NoOtherCond} Let \(Q(z)\) be a \(\mathbb C^{m\times m}\)-valued
  function,  holomorphic and single-valued in a punctured neighborhood \(\Omega\)
  of a point \(t\).
  Let
the point \(t\) be a simple pole for   \(Q(z)\), let
\begin{equation}\label{lormr1}
 Q(z)=\frac{Q_t}{z-t}+C+o(1)\text{ as }z\to t\end{equation}
  be the Laurent
expansion of the function \(Q\) at the point \(t\) and let \(R\) be a fundamental
solution  of the linear differential system
\begin{equation}\label{Fueq1}
R^\prime(z)=Q(z)R(z),\quad z\in\Omega.
\end{equation}
Assume that one of the following two cases takes place.
\begin{enumerate}
\item The coefficients \(Q_t,C\) of the expansion \eqref{lormr1} satisfy the
relations
\begin{subequations}\label{necpo01}
\begin{align}\label{necpo11}
 Q_t^2&=-Q_t,\\
 Q_tCQ_t&=-CQ_t.\label{necpo31} \end{align}
\end{subequations}

\item The coefficients \(Q_t,C\) of the expansion \eqref{lormr1} satisfy  the
relations
\begin{subequations}\label{necze01}
\begin{align}\label{necze11}
 Q_t^2&=Q_t, \\ \label{necze31}
Q_tCQ_t &=Q_tC.\end{align}
\end{subequations}
\end{enumerate}
Then \(R\) is a {\em single-valued}  function in \(\Omega\) and \(t\) is a simple
singular point of \(R\); in the first case \(t\) is a pole of \(R\), in the second
case \(t\) is a zero of \(R\).
\end{proposition}%
\begin{proof}
Once again, we shall prove only the first statement. Thus we assume that the
relations \eqref{necpo1}, \eqref{necpo3} hold and consider the transformation
\[
U(z)=(I+Q_t+(z-t)Q_t)R(z).\] Then, because of \eqref{necpo11}, the inverse
transformation is given by
 \[ R(z)=(I+Q_t+(z-t)^{-1}Q_t)U(z).\]
 Substituting these formulae
  and the Laurent expansion of
\(M\) into the linear system \eqref{Fueq1}, we obtain  the following linear system
for \(U\):
\begin{equation*}
U^{\prime}(z)= \left(\frac{(I+Q_t)CQ_t}{z-t}+V(z)\right)U(z),
\end{equation*}
where the function \(V(z)\) is  holomorphic in the entire (non-punctured)
neighborhood of the point \(t\). In view of \eqref{necpo31},
 the coefficients of this system are holomorphic at the point \(t\),
hence \(U\) is holomorphic and invertible in the entire neighborhood of \(t\) and
\(R\) has a simple pole at \(t\). Since
\[R^{-1}(z)=U^{-1}(z)(I+Q_t+(z-t)Q_t),\]
\(R^{-1}\) is holomorphic at \(t\) and hence has a zero at \(t\).
\end{proof}

 An important role in the theory of Fuchsian differential systems is played by
multiplicative decompositions of fundamental solutions (see Section 5 of
\cite{KaVo}). In the present setting we are interested in decompositions of the
following form:

\begin{definition}\label{DefPrinc}
Let \(R(z)\) be a \(\mathbb{C}^{m\times m}\)-valued function, holomorphic and
invertible in a punctured neighborhood \(\Omega\) of a point \(t\). Let \(R\) admit
in \(\Omega\) the factorization
\begin{equation}\label{muldeco}
R(z)=H_t(z)E_t(\zeta),\quad \zeta=z-t,\quad  z\in\Omega,
\end{equation}
where the factors
 \(H_t(z)\) and \(E_t(\zeta)\) possess the following properties:
 \begin{enumerate}
 \item \(H_t(z)\) is a \(\mathbb{C}^{m\times m}\)-valued function, holomorphic and
invertible  in the entire neighborhood \(\Omega\cup\{t\}\) ;

\item \(E_t(\zeta) \)
 is a \(\mathbb{C}^{m\times m}\)-valued function, holomorphic and
invertible  in the punctured plane \(\mathbb{C}_*=\mathbb{C}\setminus {0}\).
\end{enumerate}
 Then
the functions   \(E_t(\zeta) \) and  \(H_t(z)\) are said to be, respectively, the
{\em principal} and {\em regular factors} of \(R\) at \(t\).
\end{definition}

\begin{remark}\label{gauprr}
The multiplicative decomposition \eqref{muldeco}, which appears in Definition
\ref{DefPrinc}, is always possible. This follows, for example,  from the results due
to G.D.Birkhoff (see \cite{Birk1}). The principal factor \(E_t(\zeta)\) is, in a
sense, the multiplicative counterpart of the principal part of the additive
(Laurent) decomposition: it contains the information about the nature of the
singularity \(t\) of \(R\).
 Of course, the principal and regular factors at the point \(t\) are determined only up to
the transformation
\begin{equation}\label{gauprin}
E_t(\zeta)\rightarrow M(\zeta)E_t(\zeta),\quad  H_t(z)\rightarrow H_t(z)M^{-1}(z-t),
\end{equation}
 where \(M(z)\) is an  invertible entire \(\mathbb
C^{m\times m}\)-valued function. However, once the choice of the principal factor
\(E_t\) is fixed, the regular factor \(H_t\) is uniquely determined and vice-versa.
\end{remark}

A possible choice of the principal factor of the function \(R\) at its simple
singular point \(t\) is described in the following

\begin{lemma}%
\label{RegPrinFact}%
Let \(R(z)\) be a \(\mathbb{C}^{m\times m}\)-valued function,
holomorphic and invertible in a punctured neighborhood of a point
\(t\) and let \(t\) be a simple singular point of \(R\). Then a
principal factor \(E_t(\zeta)\) of \(R\) at \(t\) can be chosen as
follows.
\begin{enumerate}
\item If \(t\) is a pole of \(R\), choose any matrix  \(L\in\mathbb C^{m\times m}\),
 satisfying the conditions
\begin{equation}\label{Prpfcon}
L^2=-L ,\quad\nul(L)=\nul(R_t),\end{equation}
where \(R_t\) is the residue of \(R\) at \(t\), and set
for \(\zeta\in\mathbb{C}_*\)
\begin{equation}\label{Pprinf}
E_t(\zeta)=I+L-{\zeta}^{-1}L.\end{equation}

\item If \(t\) is a zero of \(R\), choose any matrix \(L\in\mathbb C^{m\times m}\),
satisfying the conditions
\begin{equation}\label{Zrpfcon}
L^2=L,\quad  \im(L)=\im(R_t),\end{equation} where \(R_t\) is the residue of \(R^{-1}\)
 at \(t\), and set
for \(\zeta\in\mathbb{C}_*\)
\begin{equation}\label{Zprinf}
E_t(\zeta)=I-L+\zeta L.\end{equation}
\end{enumerate}
 \end{lemma}

\begin{proof}
Let us assume that \(t\) is a pole of \(R\) and that the function \(E_t\) is given
by \eqref{Pprinf}, where the matrix \(L\) satisfies the conditions \eqref{Prpfcon}.
Then \(E_t(\zeta)\) is holomorphic in \(\mathbb C^*\); its inverse
\(E_t^{-1}(\zeta)\) is given by
\[E_t^{-1}(\zeta)=I+L-{\zeta}{L}\]
and is holomorphic in \(\mathbb C^*\), as well. Let us now show that the function
\begin{equation*}%
H(z)\stackrel{\textup{\tiny def}}{=}R(z)E_t^{-1}(z-t)
\end{equation*}%
is holomorphic and invertible at \(t\).

Indeed, in a neighborhood of \(t\) the principal part of the Laurent expansion of
\(H\) equals to \(\dfrac{R_t(I+L)}{z-t}\). But by \eqref{Prpfcon}
\(\im(L^*)=\im(R_t^*)\) and hence
\[\im((I+L^*)R_t^*)=\im((I+L^*)L^*)=\im((L^2+L)^*)=\{0\}.\]
Therefore, \(R_t(I+L)=0\) and \(H\) is holomorphic
 at \(t\).

 In the same way,  the principal part of the Laurent expansion of
\(H^{-1}\) equals to \(-\dfrac{LB_0}{z-t}\), where \(B_0=R^{-1}(t)\) is the constant
term of the Laurent expansion of \(R^{-1}\) at \(t\). But \(R_tB_0=0\) (see
\eqref{later1} in the proof of Proposition \ref{reslogder}), hence
\[\im(B_0^*L^*)=\im(B_0^*R_t^*)=\{0\},\]
\(LB_0=0\) and \(H^{-1}\) is holomorphic
 at \(t\), as well.

 The proof in the case when \(t\) is a zero of \(R\) is completely analogous.
 \end{proof}

\begin{remark}%
\label{WhyPrinc}%
Let us note  that the formulae \eqref{Pprinf} and \eqref{Zprinf} can be rewritten in
the unified form \[E_t(\zeta)=\zeta^L(=e^{L\log\zeta}).\] This is precisely the form
of the principal factor (with \(\hat Q=0\)) which appears in Proposition 5.6 of
\cite{KaVo}.
\end{remark}%
\begin{remark}\label{semi2}
The relations \eqref{Prpfcon} mean that \(-L^*\) is a projector onto \(\im(R_t^*)\).
This is equivalent to \(L\) being of the form \(L=p g_t\), where \(g_t\) is
   the right semiresidue of the function \(R\)
at its pole \(t\) and \(p\in\mathbb{C}^{m\times\rank(R_t)}\) is such that
\(g_tp=-I\). Analogously, the relations \eqref{Zrpfcon} mean that \(L\) is a
projector onto \(\im(R_t)\). This is equivalent to \(L\) being of the form
\(L=f_tq\), where \(f_t\) is
   the left semiresidue of the function \(R^{-1}\)
at its pole \(t\) and \(q\in\mathbb{C}^{\rank(R_t)\times m}\) is such that
\(qf_t=I\). For example, one can choose the matrix \(L\) mentioned in Lemma
\ref{RegPrinFact} as follows:
\[L=\left\{\begin{array}{l@{\quad\text{if}\quad}l}
-g_t^*(g_tg_t^*)^{-1}g_t& t\text{ is a pole of }R,\\[1ex]
\phantom{-}f_t(f_t^*f_t)^{-1}f_t^*& t\text{ is a zero of }R.
\end{array}\right.
\]
\end{remark}

\section{Rational matrix functions of simple structure}
In this section we apply the local results obtained in Section \ref{GSMMF} to the
study of rational matrix functions.
\begin{definition}\label{simstruc}
 A \(\mathbb C^{m\times m}\)-valued
 rational function \(R(z)\) is said to be a rational matrix function
{\em of simple structure} if  it meets  the following conditions:
\begin{enumerate}
\item
 \(\det
R(z)\not\equiv 0\); \item   all singular points of \(R\) are simple; \item
\(z=\infty\) is a regular point of \(R\).
\end{enumerate}
 The
set of all poles of the function \(R\) is said to be the {\em pole set} of the
function \(R\) and is denoted by \(\mathcal{P}_{R}\). The set of all zeros of the function
\(R\) is said to be the {\em zero set} of the function \(R\) and is denoted by
\(\mathcal{Z}_{R}\).
\end{definition}

\begin{remark}
Note that  if \(R\) is a rational matrix function  of simple structure then the
inverse function \(R^{-1}\) is a rational matrix function of simple structure, as
well, and \(\mathcal{Z}_{R}=\mathcal{P}_{R^{-1}}\).
\end{remark}

 Below we
formulate the "global" counterparts of Propositions \ref{reslogder} and
\ref{NoOtherCond} in order to  characterize Fuchsian differential systems whose
fundamental solutions are rational matrix functions of simple structure.

\begin{theorem}%
\label{CharLogDer}%
Let \(R(z)\) be a rational matrix function of simple structure with the pole set
\(\mathcal{P}_{R}\) and the zero set \(\mathcal{Z}_{R}\). Then its logarithmic
derivative\footnote{See Remark \ref{leftder}.} \(Q_R(z)\)  is a rational function
with the set of poles \(\mathcal{P}_{R}\cup\mathcal{Z}_{R}\);  all the poles of \(Q_R\)
are simple. Furthermore, the function \(Q_R\) admits the additive decomposition
\begin{equation}
\label{AELD}%
Q_R(z)=\sum_{t\in\mathcal{P}_{R}\cup\mathcal{Z}_{R}}\frac{Q_t}{z-t},
\end{equation}
and its residues \(Q_t\in\mathbb C^{m\times m}\)  satisfy the following
relations:
\begin{equation}
\label{rlogd1} \sum_{t\in\mathcal{P}_{R}\cup\mathcal{Z}_{R}}Q_t=0,\end{equation}
\begin{align} \label{rlogd2}
Q_t^2&=\left\{\begin{array}{l@{\quad\text{if}\quad}l} -Q_t & t\in\mathcal{P}_{R},\\
\phantom{-}Q_t & t\in\mathcal{Z}_{R},\end{array}\right.\\
\label{rlogd3}
Q_tC_tQ_t&=\left\{\begin{array}{l@{\quad\text{if}\quad}l} -C_tQ_t & t\in\mathcal{P}_{R},\\
\phantom{-}Q_tC_t & t\in\mathcal{Z}_{R},\end{array}\right.
\end{align}
where
\begin{equation}\label{conterq}
C_{t}=
\sum_{\substack{t^{\prime}\in\mathcal{P}_{R}\cup\mathcal{Z}_{R}\\
t^{\prime}\not=t}}%
\frac{Q_{t^\prime}}{t-t^\prime}.
\end{equation}%
\end{theorem}

\begin{proof}
Since both functions \(R\) and \(R^{-1}\) are holomorphic
 in \(\overline{\mathbb{C}}\setminus(\mathcal{P}_{R}\cup\mathcal{Z}_{R}),\)
 the logarithmic derivative \(Q_R\) is holomorphic there, as well.
According to Proposition \ref{reslogder}, each point of the set
\(\{\mathcal{P}_{R}\cup\mathcal{Z}_{R}\}\) is a simple pole of \(Q_R\), hence  we can
write for \(Q_R\) the additive decomposition
\[Q_R(z)=Q_R(\infty)+\sum_{t\in\mathcal{P}_{R}\cup\mathcal{Z}_{R}}\frac{Q_t}{z-t},\]
 where \(Q_t\) are the residues of \(Q_R\).
 Since \(R\) is holomorphic at \(\infty\), the entries of its derivative \(R^\prime\)
 decay as \(o(|z|^{-1})\) when \(z\rightarrow\infty\).
 The rate of decay for the logarithmic derivative \(Q_R\) is the same, because
  \(R^{-1}\), too, is  holomorphic at \(\infty\).
  Thus we obtain  the additive decomposition \eqref{AELD} for \(Q_R\)
  and the relation \eqref{rlogd1} for the residues \(Q_t\).
  Now the relations \eqref{rlogd2}, \eqref{rlogd3}
  follow immediately from Proposition \ref{reslogder}, once we observe that
  the matrix \(C_t\) given by \eqref{conterq} is but the constant term
   of the Laurent expansion of \(Q_R\) at its pole \(t\).
\end{proof}

\begin{theorem} %
\label{CharLogDer2}%
Let \(\mathcal{P}\) and \(\mathcal{Z}\) be two finite  disjoint subsets of
the complex plane \(\mathbb{C}\) and  let \(Q(z)\) be a \(\mathbb{C}^{m\times
m}\)-valued rational function of the form
\begin{equation}%
\label{DELD} Q(z)=\sum_{t\in\mathcal{P}\cup\mathcal{Z}}\frac{Q_t}{z-t},
\end{equation}
where  \(Q_t\in\mathbb C^{m\times m}.\)  Let the matrices \(Q_t\) satisfy the
relations
\begin{equation}
\label{rlogd11} \sum_{t\in\mathcal{P}\cup\mathcal{Z}}Q_t=0,\end{equation}
\begin{align}
\label{rlogd21}
Q_t^2&=\left\{\begin{array}{l@{\quad\text{if}\quad}l} -Q_t & t\in\mathcal{P},\\
\phantom{-}Q_t & t\in\mathcal{Z},\end{array}\right.\\
\label{rlogd31}
Q_tC_tQ_t&=\left\{\begin{array}{l@{\quad\text{if}\quad}l} -C_tQ_t & t\in\mathcal{P},\\
\phantom{-}Q_tC_t & t\in\mathcal{Z},\end{array}\right.
\end{align}
where
\begin{equation}\label{conterq1}
C_{t}=
\sum_{\substack{t^{\prime}\in\mathcal{P}\cup\mathcal{Z}\\
t^{\prime}\not=t}}%
\frac{Q_{t^\prime}}{t-t^\prime}.
\end{equation}%
 Let \(R(z)\) be a fundamental solution of the
Fuchsian differential system
\begin{equation}\label{Fueqg}
R^\prime(z)=Q(z)R(z).
\end{equation}
Then \(R\) is  a rational matrix function  of simple structure such that
\[\mathcal{P}_{R}=\mathcal{P},\quad \mathcal{Z}_{R}= \mathcal{Z}.\]
\end{theorem}

\begin{proof}
Since the condition \eqref{rlogd11} implies that the point \(\infty\) is a regular
point for the Fuchsian system \eqref{Fueqg}, we may, without loss of generality,
consider the fundamental solution \(R\) satisfying the initial condition
\(R(\infty)=I\). Then  \(R(z)\) is a matrix function, holomorphic (a priori,
multi-valued) and invertible in the (multi-connected) set
\(\overline{\mathbb{C}}\setminus(\mathcal{P}\cup\mathcal{N})\). However, for
\(t\in\mathcal{P}\cup\mathcal{N}\) the function \(Q\) admits in a neighborhood of
\(t\) the Laurent expansion
\[Q(z)=
\frac{Q_t}{z-t}+C_t+ o(1)\] with the constant term \(C_t\) given by
\eqref{conterq1}. The coefficients \(Q_t\) and \(C_t\) satisfy the relations
\eqref{rlogd21}, \eqref{rlogd31}, hence by Proposition \ref{NoOtherCond} the
function \(R\) is meromorphic at \(t\). Since this is true for every
\(t\in\mathcal{P}\cup\mathcal{N}\), the function \(R\) is rational (in particular,
single-valued).  Proposition \ref{NoOtherCond} also implies that every
\(t\in\mathcal{P}\) (respectively, \(t\in\mathcal{Z}\)) is a simple pole
(respectively, a zero) of the function \(R\) and a zero (respectively, a simple
pole) of the inverse function \(R^{-1}\). Therefore, \(R\) is  a rational matrix
function  of simple structure with the pole set \(\mathcal{P}\) and the zero set
\(\mathcal{Z}\).
\end{proof}

We close this section with the following useful

\begin{lemma}%
\label{DimRFSSLem}%
Let \(R\) be a rational matrix function of simple structure. For
\(t\in\mathcal{P}_{R}\cup\mathcal{Z}_{R}\) let \(R_t\) denote the residue of the
function \(R\) at  \(t\) if  \(\,t\in\mathcal{P}_{R}\), and the residue of the
inverse function \(R^{-1}\) at \(t\) if  \(\,t\in\mathcal{Z}_{R}\). Then
\begin{equation}
\label{DimRFSS}
\sum_{t\in\mathcal{P}_{R}}\rank(R_t)=
\sum_{t\in\mathcal{Z}_{R}}\rank(R_t).
\end{equation}
\end{lemma}%

\begin{proof}
Let us consider the logarithmic derivative \(Q_R\) of \(R\). Its residues \(Q_t\)
satisfy by Theorem \ref{CharLogDer} the relations \eqref{rlogd1} and \eqref{rlogd2}.
From \eqref{rlogd2} it follows that
\[\rank(Q_t)=\left\{\begin{array}{l@{\quad\text{if}\quad}l}
-\trace(Q_t)& t\in\mathcal{P}_{R},\\
\phantom{-}\trace(Q_t)& t\in\mathcal{Z}_{R}.
\end{array}\right.\]
But \eqref{rlogd1} implies
\[\sum_{t\in\mathcal{P}_{R}}\trace(Q_t)+
\sum_{t\in\mathcal{Z}_{R}}\trace(Q_t)=0,\]
hence
\[\sum_{t\in\mathcal{P}_{R}}\rank(Q_t)=
\sum_{t\in\mathcal{Z}_{R}}\rank(Q_t).\]
Finally, by Proposition \ref{reslogder} (see \eqref{necpo2}, \eqref{necze2} there),
\[\rank(R_t)=\rank(Q_t),\quad \forall t\in\mathcal{P}_{R}\cup\mathcal{Z}_{R}.\]
Thus \eqref{DimRFSS} holds.
\end{proof}

\section{Generic rational matrix functions}
\label{RMFGP}
\begin{definition}%
\label{RatFunGenRosDef} A \(\mathbb{C}^{m\times m}\)-valued rational function \(R(z)\)
  is said to be a {\em generic}\footnote{ In \cite{Kats2} such functions are called
  "rational matrix functions in general position".}
rational matrix function  if \(R\) is a rational matrix function of simple structure
 and  all the residues  of the functions \(R\) and \(R^{-1}\)
have rank one.
\end{definition} %

\begin{lemma}%
\label{CardCoLem}%
Let \(R\) be a generic rational matrix function.
 Then the cardinalities of its pole
and zero sets\footnote{ See Definition \ref{simstruc}.} coincide:
\begin{equation}%
\label{CardCo}%
 \#\mathcal{P}_{R}=\#\mathcal{Z}_{R}.
\end{equation}%
\end{lemma}  %

\begin{proof} Since all the residues  of  \(R\) and \(R^{-1}\)
are of rank one, the statement follows immediately
from Lemma \ref{DimRFSSLem}.\end{proof}

Let \(R\) be a \(\mathbb{C}^{m\times m}\)-valued generic rational function. In what
follows, we assume that \(R\) is normalized by
\begin{equation}\label{norm}
R(\infty)=I.
\end{equation}
Let us order  somehow the pole and zero sets of \(R\):
 \begin{equation}\label{tpzs}
\mathcal{P}_{R}=\{t_1, \dots, t_n \},\quad
 \mathcal{Z}_{R}=\{t_{n+1}, \dots, t_{2n} \},\end{equation}
where \(n=\#\mathcal{P}_{R}=\#\mathcal{Z}_{R}\).
Then we can write for \(R\) and \(R^{-1}\) the
additive decompositions
\begin{subequations}
\label{aer}
\begin{align}
\label{AdExp}
R(z)&=I+\sum_{k=1}^n\frac{R_k}{z-t_k},\\
R^{-1}(z)&=I+\sum_{k=n+1}^{2n}\frac{R_k}{z-t_k},
\label{AdExpInv}
\end{align}
\end{subequations}
where for \(1\leq k\leq n\) (respectively, \(n+1\leq k\leq 2n\)) we denote by
\(R_k\) the residue of \(R\) (respectively, \(R^{-1}\)) at its pole \(t_k\). Since
each matrix \(R_k\) is of rank one,
 the representations \eqref{aer}
  can be rewritten as
\begin{subequations}
\label{afr}
\begin{align}
R(z)&=
I+\sum_{k=1}^nf_{k}\frac{1}{z-t_k}
g_k,
\label{AddFactDir}\\
R^{-1}(z)&=I+\sum_{k=n+1}^{2n}f_k\frac{1}{z-t_k} g_k,
\label{AddFactInv}%
\end{align}
\end{subequations}
 where for \(1\leq k \leq n\) (respectively, \(n+1\leq k \leq 2n\))
  \(f_k\in\mathbb{C}^{m\times 1}\) and
 \(g_k\in\mathbb{C}^{1\times m}\) are the left and right semiresidues
 \footnote{ See Definition \ref{DefSimPol}.} of \(R\)
 (respectively, \(R^{-1}\)) at \(t_k\).
Furthermore,
 we introduce two \(n\times n\) diagonal matrices:
\begin{equation}
A_{\mathcal{P}}=\diag(t_1, \dots ,t_n),
\quad A_{\mathcal{Z}}=\diag(t_{n+1}, \dots, t_{2n}),
\label{polematr}
\end{equation}
 two \( m\times n \) matrices :
\begin{equation}
F_{\mathcal{P}} = \begin{pmatrix}f_ 1 & \dots & f_ n\end{pmatrix} ,
\quad F_{\mathcal{Z}} = \begin{pmatrix}f_{n+ 1} & \dots & f_{2n}\end{pmatrix} ,
\label{leftsemires}
\end{equation}
and two \( n\times m \) matrices :
\begin{equation}
G_{\mathcal{P}} =\begin{pmatrix} g_ 1 \\ \vdots \\ g_ n\end{pmatrix} ,
\quad G_{\mathcal{Z}} = \begin{pmatrix}g_{n+ 1} \\ \vdots \\ g_{2n}\end{pmatrix} .
 \label{rightsemires}
\end{equation}
The matrices \( A_{\mathcal{P}}\) and \( A_{\mathcal{Z}}\) are said to be,
respectively, the {\em pole} and {\em zero matrices} of \(R \). The matrices \(
F_{\mathcal{P}} \) and \( G_{\mathcal{P}} \) are said to be, respectively, the {\em
left} and {\em right pole semiresidual matrices} of \(R \). Analogously, the
matrices \( F_{\mathcal{Z}}\) and \(G_{\mathcal{Z}} \) are said to be the {\em left}
and {\em right zero semiresidual matrices} of \(R \).

\begin{remark}\label{gauge}
It should be mentioned that for a fixed ordering \eqref{tpzs}
of the pole and zero sets
 the pole and the zero matrices \(A_{\mathcal{P}}\) and
\(A_{\mathcal{N}}\) are defined uniquely, and the semiresidual matrices
\(F_{\mathcal{P}}, G_{\mathcal{P}} , F_{\mathcal{Z}}, G_{\mathcal{Z}}\)
 are defined essentially uniquely, up to the transformation
\begin{subequations}
\begin{align} F_{\mathcal{P}}\rightarrow F_{\mathcal{P}} D_{\mathcal{P}},
&\quad G_{\mathcal{P}}\rightarrow D_{\mathcal{P}}^{-1} G_{\mathcal{P}},
\label{polegauge}
\\
F_{\mathcal{Z}}\rightarrow F_{\mathcal{Z}}  D_{\mathcal{Z}}, &\quad
G_{\mathcal{Z}}\rightarrow D_{\mathcal{Z}}^{-1} G_{\mathcal{Z}}, \label{zerogauge}
 \end{align}
 \end{subequations}
where \(D_{\mathcal{P}},D_{\mathcal{Z}}\in\mathbb{C}^{n\times n}\) are arbitrary
invertible diagonal matrices. Once the choice of the left pole semiresidual matrix
\( F_{\mathcal{P}}\) is fixed, the right pole semiresidual matrix
\(G_{\mathcal{P}}\) is determined uniquely, etc.
\end{remark}

In terms of the matrices
\(A_{\mathcal{P}},\,A_{\mathcal{Z}},\,F_{\mathcal{P}},G_{\mathcal{Z}},
F_{\mathcal{Z}},G_{\mathcal{Z}} \), the representations \eqref{afr} take the
following form:
\begin{subequations}
\label{aemf}
\begin{align}%
\label{AddDirMaF}%
R(z)&=I+F_{\mathcal{P}}\left(zI-A_{\mathcal{P}}\right)^{-1}G_{\mathcal{P}},\\
\label{AddInvMaF}%
R^{-1}(z)&=I+F_{\mathcal{Z}}\left(zI-A_{\mathcal{Z}}\right)^{-1} G_{\mathcal{Z}}.
\end{align}%
\end{subequations}

The representations \eqref{aemf} are not quite satisfactory for the
following reasons. Firstly, in view of the identity
\(RR^{-1}=R^{-1}R=I,\) the information contained in the pair of
representations \eqref{aemf} is {\em redundant}: each of these
representations determines the function \(R\) (and \(R^{-1}\))
uniquely. Secondly, if, for example, the diagonal matrix
\(A_{\mathcal{P}}\) and the matrices
\(F_{\mathcal{P}},G_{\mathcal{P}}\) of appropriate dimensions are
chosen arbitrarily then the rational function \eqref{AddDirMaF} need
not be generic. In our investigation we shall mainly use another
version of the system representation of rational matrix functions,
more suitable for application to linear differential equations. This
is the so-called joint system representation (see \cite{Kats2} for
details and references) of the function \(R(z)R^{-1}(\omega)\)  of
{\em two} independent variables \(z\) and \(w\). A key role in the
theory of the  joint system representation is played by
 the {\em Lyapunov equations}. These are matricial
equations of the form
\begin{equation}%
\label{LyEq}%
UX-XV=Y,
\end{equation}%
where the matrices \(U,V,Y\in\mathbb{C}^{n\times n}\) are given, and the matrix
\(X\in\mathbb{C}^{n\times n}\) is unknown. If the spectra of the matrices \(U\) and
\(V\) are disjoint, then the Lyapunov equation \eqref{LyEq} is uniquely solvable
with respect to \(X\) for arbitrary right-hand side \(Y\). The solution \(X\) can be
expressed, for example, as the contour integral
\begin{equation}
\label{CoInR}%
 X=\frac{1}{2\pi i}\oint\limits_{\Gamma}(z
I-U)^{-1}Y(z I-V)^{-1}dz,
\end{equation}
where \(\Gamma\) is an arbitrary contour, such that the spectrum of
\(U\) is inside \(\Gamma\) and the spectrum of \(V\) is outside
\(\Gamma\) (see, for instance, Chapter I, Section 3 of the book
\cite{DaKr}).

With the generic rational function \(R\) we associate the following pair of Lyapunov
equations (with unknown \(S_{\mathcal{ZP}},S_{\mathcal{PZ}}\in\mathbb{C}^{n\times
n}\)):
\begin{subequations}\label{SL}
\begin{align}
A_{\mathcal{Z}}S_{\mathcal{ZP}}-S_{\mathcal{ZP}}A_{\mathcal{P}}&=G_{\mathcal
Z}F_{\mathcal{P}}, \label{SL3}\\
A_{\mathcal{P}}S_{\mathcal{PZ}}-S_{\mathcal{PZ}}A_{\mathcal{Z}}&=G_{\mathcal
{P}}F_{\mathcal{Z}}.
\label{SL4}%
\end{align}
\end{subequations}
Since the spectra of the pole and zero matrices \(A_{\mathcal{P}}\) and
\(A_{\mathcal{Z}}\) do not intersect (these  are  the pole  and zero sets of \(R\)),
the Lyapunov equations \eqref{SL3} and \eqref{SL4} are uniquely solvable. In fact,
since the matrices \(A_{\mathcal{P}}\) and \(A_{\mathcal{Z}}\) are diagonal, the
solutions can be given explicitly (using the notation  \eqref{polematr} --
\eqref{rightsemires}):
\begin{equation}%
\label{ESNP}%
S_{\mathcal{ZP}}=\left(\frac{g_{n+i}f_{j}}{t_{n+i}-t_j} \right)_{i,j=1}^n,\quad
S_{\mathcal{PZ}}=\left(\frac{g_{i}f_{n+j}}{t_i-t_{n+j}} \right)_{i,j=1}^n.
\end{equation}%
The matrices \(S_{\mathcal{ZP}}\) and \(S_{\mathcal{PZ}}\) are said to be,
respectively, the {\em zero-pole} and {\em pole-zero coupling matrices} of
  \(R\).

\begin{proposition} %
\label{MuInc}%
Let \(R(z)\) be a generic rational matrix function normalized by \(R(\infty)=I.\)
Then
\begin{enumerate}
\item
 the coupling matrices
\(S_{\mathcal{ZP}}\) and \(S_{\mathcal{PZ}}\) of  \(R\) are mutually inverse:
\begin{equation}%
\label{MIR}%
S_{\mathcal{ZP}} S_{\mathcal{PZ}}=S_{\mathcal{PZ}} S_{\mathcal{ZP}}=I;
\end{equation}%
\item  for the semiresidual matrices of \(R\) the following relations hold:
\begin{equation}%
\label{ZPCR}%
G_{\mathcal{Z}}=-S_{\mathcal{ZP}}G_{\mathcal{P}},\quad F_{\mathcal{P}}=F_{\mathcal
Z}S_{\mathcal{ZP}};
\end{equation}
\item the function \(R\) admits the joint representation
\begin{equation}\label{jointrep0}
R(z)R^{-1}(\omega)=I+(z-\omega)F_{\mathcal{P}}
(zI-A_{\mathcal{P}})^{-1}S_{\mathcal{ZP}}^{-1} (\omega
I-A_{\mathcal{Z}})^{-1}G_{\mathcal{Z}},
\end{equation}%
where \(A_{\mathcal{P}},A_{\mathcal{Z}}\) are the pole and zero matrices of \(R\).
\end{enumerate}
\end{proposition}%
\begin{proof}
The proof of Proposition \ref{MuInc} can be found in \cite{Kats2}.
\end{proof}
\begin{remark}\label{separate}
Note that, since \(R(\infty)=I\), one can recover from the joint representation
\eqref{jointrep0} when \(z\to\infty\) or \(\omega\to\infty\) the separate
representations
\begin{subequations}\label{amf1}
\begin{align}
\label{admf1} R(z)&=I-F_{\mathcal{P}} (zI-A_{\mathcal{P}})^{-1}S_{\mathcal{ZP}}^{-1}
G_{\mathcal{Z}},\\
\label{aimf1} R^{-1}(\omega)&=I+F_{\mathcal{P}} S_{\mathcal{ZP}}^{-1} (\omega
I-A_{\mathcal{Z}})^{-1}G_{\mathcal{Z}},
\end{align}
\end{subequations}
which, in view of \eqref{ZPCR}, coincide with the representations \eqref{aemf}.
\end{remark}

\begin{remark}\label{secjoint}
In view of \eqref{MIR}, \eqref{ZPCR}, one can also write the joint representation
for \(R\) in terms of the matrices \(F_{\mathcal{Z}},G_{\mathcal{P}}\) and the
solution \(S_{\mathcal{PZ}}\) of the Lyapunov equation \eqref{SL4}:
\begin{equation}\label{jointrep1}
R(z)R^{-1}(\omega)=I-(z-\omega)F_{\mathcal{Z}}S_{\mathcal{PZ}}^{-1}
(zI-A_{\mathcal{P}})^{-1}S_{\mathcal{PZ}} (\omega
I-A_{\mathcal{Z}})^{-1}S_{\mathcal{PZ}}^{-1}G_{\mathcal{P}}.
\end{equation}%
Thus we may conclude that the pole and zero sets together with a pair of the
semiresidual matrices (either right pole and left zero or left pole and right zero)
determine the normalized generic rational function \(R\) {\em uniquely}.
\end{remark}

\begin{remark}
 The theory of system representations for rational matrix functions with
prescribed zero and pole structures first appeared in \cite{GKLR}, and was further
developed in \cite{BGR1}, \cite{BGR2}, and \cite{BGRa}.

The joint representations \eqref{jointrep0}, \eqref{jointrep1} suggest that this
theory can be  applied to the investigation of families
 of rational
functions parameterized by the zeros' and poles' loci and the
corresponding\footnote{ See Propositions \ref{CharLogDer}, \ref{CharLogDer2}.}
deformations of linear differential systems.
 The version of
this theory adapted for such applications was presented in \cite{Kats1} and
\cite{Kats2}. Also in \cite{Kats2} one can find some historical remarks and a list
of references.
\end{remark}

\begin{proposition}
\label{ReprLogDer} Let \(R(z)\) be a generic rational matrix function normalized by
\(R(\infty)=I.\) Then its logarithmic derivative\footnote{See Remark \ref{leftder}.}
  admits the representation:
\begin{equation}\label{replogder}
R^{\prime}(z)R^{-1}(z)= F_{\mathcal{P}} (zI-A_{\mathcal{P}})^{-1}S_{\mathcal{ZP}}^{-1}
(zI-A_{\mathcal{Z}})^{-1}G_{\mathcal{Z}}
\end{equation}
where \(A_{\mathcal{P}}\) and \(A_{\mathcal{Z}}\) are the pole and zero matrices of
\(R\); \(F_{\mathcal{P}}\) and \(G_{\mathcal{Z}}\) are the left pole and right zero
semiresidual matrices of \(R\); \(S_{\mathcal{ZP}}\) is the zero-pole coupling
matrix of \(R\).
\end{proposition}%

\begin{proof}
Differentiating \eqref{jointrep0} with respect to \(z\), we obtain
\[
R^\prime(z)R^{-1}(\omega)=F_{\mathcal{P}}\left(I-(z-\omega)
(zI-A_{\mathcal{P}})^{-1}\right) (zI-A_{\mathcal{P}})^{-1}S_{\mathcal{ZP}}^{-1}
(\omega I-A_{\mathcal{Z}})^{-1}G_{\mathcal{Z}}. \] Now set \(\omega=z\) to obtain
\eqref{replogder}.
\end{proof}

\begin{remark}\label{seclogder}
The representation \eqref{replogder} for the logarithmic derivative \(Q_R\) of the
normalized generic rational matrix function \(R\)   can also be rewritten in terms
of the matrices \(F_{\mathcal{Z}},G_{\mathcal{P}}\) and the solution
\(S_{\mathcal{PZ}}\) of the Lyapunov equation \eqref{SL4} (see Remark
\ref{secjoint}):
\begin{equation}\label{replogder1}
R^\prime(z)R^{-1}(z)=-F_{\mathcal{Z}}S_{\mathcal{PZ}}^{-1}
(zI-A_{\mathcal{P}})^{-1}S_{\mathcal{PZ}} (z
I-A_{\mathcal{Z}})^{-1}S_{\mathcal{PZ}}^{-1}G_{\mathcal{P}}.
\end{equation}%
\end{remark}

\section{Generic
rational matrix functions  with prescribed local data} \label{RFPD}
 In the previous
section we discussed the  question, how to represent
 a generic rational
function \(R\) in terms of its local data (the pole and zero sets and the residues).
The main goal of this section is to construct
 a  (normalized) generic rational matrix \(R(z)\)
function with {\em prescribed} local data. In view of Proposition \ref{MuInc},
Remark \ref{secjoint} and Remark \ref{gauge}, such data should be given in
 the form  of two diagonal
matrices of the same dimension (the pole and zero matrices) and two semiresidual
matrices of appropriate dimensions (either right pole and left zero or right zero
and left pole)\footnote{Here we use the terminology introduced in Section
\ref{RMFGP}.}.

Thus we consider the following

\begin{Pb}\label{pres1}
Let two diagonal matrices
\[
A_{\mathcal{P}}=\diag(t_1, \dots ,t_n), \quad A_{\mathcal{Z}}=\diag(t_{n+1}, \dots,
t_{2n}),\quad t_i\not=t_j\text{ unless } i=j,\] and two matrices:
\(F\in\mathbb{C}^{m\times n}\), \(G\in\mathbb{C}^{n\times m}\), be given.
\begin{paragraph}{\(\mathcal{ZP}\)-version:} Find a generic \(\mathbb{C}^{m\times m}\)-valued rational
function \(R(z)\) such that
\begin{enumerate}
\item \(R(\infty)=I\);

 \item the matrices \(A_{\mathcal{P}}\) and \(A_{\mathcal{Z}}\)
are, respectively,
 the
pole and zero matrices of \(R\);

\item the matrix \(F\) is the left pole semiresidual matrix \(F_\mathcal{P}\) of
\(R\): \(\,F=F_\mathcal{P}\);

\item the matrix \(G\) is the right zero semiresidual matrix \(G_\mathcal{Z}\) of
\(R\): \(\,G=G_\mathcal{Z}\).
\end{enumerate}
\end{paragraph}
\begin{paragraph}{\(\mathcal{PZ}\)-version:} Find a generic \(\mathbb{C}^{m\times m}\)-valued
rational function \(R(z)\)
 such that

\begin{enumerate}
\item \(R(\infty)=I\);

\item the matrices \(A_{\mathcal{P}}\) and \(A_{\mathcal{Z}}\) are, respectively,
 the
pole and zero matrices of \(R\);

\item the matrix \(F\) is the left zero semiresidual matrix \(F_\mathcal{Z}\) of
\(R\): \(\,F=F_\mathcal{Z}\);

\item the matrix \(G\) is the right pole semiresidual matrix \(G_\mathcal{P}\) of
\(R\): \(\,G=G_\mathcal{P}\).
\end{enumerate}
\end{paragraph}
\end{Pb}

\begin{proposition}\label{Lisnecas}\mbox{}
\begin{enumerate}
\item The \(\mathcal{ZP}\)-version of Problem \ref{pres1} is solvable if and only
if the solution \(S\) of the Lyapunov equation
\begin{equation}\label{ZPLE}
A_\mathcal{Z}S-SA_\mathcal{P}=GF
\end{equation}
 is an invertible matrix.

 \item The \(\mathcal{PZ}\)-version of Problem
\ref{pres1} is solvable if and only if the solution \(S\) of the Lyapunov equation
\begin{equation}\label{PZLE}
A_\mathcal{P}S-SA_\mathcal{Z}=GF
\end{equation}
 is an invertible matrix.
\end{enumerate}
\end{proposition}
\begin{proof} The proof of Proposition \ref{Lisnecas} can be found in \cite{Kats2}.
Here we would like to note that the  solutions of the  Lyapunov equations
\eqref{ZPLE} and \eqref{PZLE} can be written explicitly as, respectively,
\begin{equation}%
\label{ESNP1}%
S=\left(\frac{g_{i}f_{j}}{t_{n+i}-t_j} \right)_{i,j=1}^n \quad
\text{and}\quad
S=\left(\frac{g_{i}f_{j}}{t_i-t_{n+j}} \right)_{i,j=1}^n,
\end{equation}%
where \(g_i\) is the \(i\)-th row of the matrix \(G\) and \(f_j\) is the \(j\)-th
column of the matrix \(F\).
Note also that the necessity of \(S\) being invertible in both cases follows from
Proposition \ref{MuInc}.
\end{proof}

 In view of Proposition \ref{Lisnecas}, we propose the following terminology:
\begin{definition}%
\label{DefLocData}%
Let \(A_{\mathcal{P}},A_{\mathcal{Z}},F,G\) be the given data of Problem
\ref{pres1}. Then:
\begin{enumerate}
\item  the solution \(S\) of the Lyapunov equation
\begin{equation}\label{ZPLE1} A_\mathcal{Z}S-SA_\mathcal{P}=GF
\end{equation} is said to be the \(\mathcal{ZP}\)-{\em coupling} matrix related to the
data \(A_{\mathcal{P}},A_{\mathcal{Z}},F,G\);

\item the solution \(S\) of the Lyapunov equation
\begin{equation}\label{PZLE1}
A_\mathcal{P}S-SA_\mathcal{Z}=GF
\end{equation}
 is said to be the \(\mathcal{PZ}\)-{\em coupling} matrix related to the
data \(A_{\mathcal{P}},A_{\mathcal{Z}},F,G\);

 \item the data \(A_{\mathcal{P}},A_{\mathcal{Z}},F,G\) are said to be
 \(\mathcal{ZP}\)-{\em
admissible} if the \(\mathcal{ZP}\)-{coupling} matrix related to this data is
invertible;

\item  the data \(A_{\mathcal{P}},A_{\mathcal{Z}},F,G\) are said to be
\(\mathcal{PZ}\)-{\em admissible} if the \(\mathcal{PZ}\)-{coupling} matrix
related to this data is invertible.
\end{enumerate}
\end{definition}

\begin{proposition}\label{record}
Let  \(A_{\mathcal{P}},A_{\mathcal{Z}},F,G\) be the given data of Problem
\ref{pres1}.
\begin{enumerate}
\item If the data \(A_{\mathcal{P}},A_{\mathcal{Z}},F,G\) are
\(\mathcal{ZP}\)-{admissible} then the \(\mathcal{ZP}\)-version of Problem
\ref{pres1} has the unique solution \(R(z)\) given by
\begin{equation}\label{jointrep10}
R(z)=I-F(
I-A_{\mathcal{P}})^{-1} S^{-1} G,\end{equation}
where \(S\) is the \(\mathcal{ZP}\)-coupling matrix related to
the data \(A_{\mathcal{P}},A_{\mathcal{Z}},F,G\).
The logarithmic derivative of \(R\) is given by
\begin{equation}\label{replogder10}
R^{\prime}(z)R^{-1}(z)= F (zI-A_{\mathcal{P}})^{-1}S^{-1} (zI-A_{\mathcal{Z}})^{-1}G.
\end{equation}

\item If the data \(A_{\mathcal{P}},A_{\mathcal{Z}},F,G\) are
\(\mathcal{PZ}\)-{admissible} then the \(\mathcal{PZ}\)-version of Problem
\ref{pres1} has the unique solution \(R(z)\) given by
\begin{equation}\label{jointrep20}
R(z)=I+F S^{-1}(zI-A_{\mathcal{P}})^{-1}G,\end{equation}
where \(S\) is the \(\mathcal{PZ}\)-coupling matrix related to
the data \(A_{\mathcal{P}},A_{\mathcal{Z}},F,G\).
The logarithmic derivative of \(R\) is given by
\begin{equation}
\label{replogder20}
R^{\prime}(z)R^{-1}(z)= -F S^{-1}(zI-A_{\mathcal{P}})^{-1}S
(zI-A_{\mathcal{Z}})^{-1}S^{-1}G.
\end{equation}
\end{enumerate}
\end{proposition}
\begin{proof}
If the data \(A_{\mathcal{P}},A_{\mathcal{Z}},F,G\) are
\(\mathcal{ZP}\)-{admissible} then, according to Proposition \ref{Lisnecas}, the
\(\mathcal{ZP}\)-version of Problem \ref{pres1} has a solution \(R\). Then the
\(\mathcal{ZP}\)-{\em coupling} matrix \(S\) related to the data
\(A_{\mathcal{P}},A_{\mathcal{Z}},F,G\) is  the zero-pole coupling matrix of \(R\).
According to  Proposition \ref{MuInc}, the function \(R\)  admits the representation
\eqref{jointrep10} and hence\footnote{See Remark \ref{separate}.} is determined
uniquely. Now the representation \eqref{replogder10} for the logarithmic derivative
 follows from Proposition \ref{ReprLogDer}.

Analogous considerations hold also in the case when the data
\(A_{\mathcal{P}},A_{\mathcal{Z}},F,G\) are \(\mathcal{PZ}\)-{admissible} (see
Remarks \ref{secjoint}, \ref{seclogder}).
\end{proof}

It was already mentioned (see Remark \ref{gauge}) that in the definition of the
semiresidual matrices there is a certain freedom. Accordingly, certain equivalency
classes rather than individual matrices
 \(F,G\) should serve as data for Problem \ref{pres1}. The
appropriate definitions are similar to the definition of the complex projective
space \(\mathbb{P}^{k-1}\) as the space of equivalency classes
 of the set
\(\mathbb{C}^{k}\setminus\{0\}\) (two vectors
\(h^{\prime},\,h^{\prime\prime}\in\mathbb{C}^{k}\setminus\{0\}\) are declared to be
{\em equivalent} if \(h^{\prime},\,h^{\prime\prime}\) are {\em proportional}, i.e.
\(h^{\prime\prime}=d h^{\prime}\) for some \(d\in\mathbb{C}_*\)).

\begin{definition}%
\label{ProjEquivDef}%
\begin{enumerate}%
\item Let \(\mathbb{C}^{m\times n}_{\ast,c}\) denote the set of \(m\times n\)
matrices which have no zero columns. Two matrices
\(F^{\prime},F^{\prime\prime}\in\mathbb{C}^{m\times n}_{\ast,c}\)
 are declared to be equivalent:
 \(F^{\prime}\stackrel{c}{\sim}
F^{\prime\prime}\), if%
\begin{equation}
\label{Eqc}%
 F^{\prime\prime}=F^{\prime}D_c,
\end{equation}
where \(D_c\) is a diagonal invertible matrix.\\
The space \(\mathbb{P}_{c}^{(m-1)\times n}\) is a factor-set of the set
\(\mathbb{C}^{m\times n}_{\ast,c}\) modulo the equivalency relation
\(\stackrel{c}{\sim}\).

\item Let \(\mathbb{C}^{n\times m}_{\ast,r}\) denote the set of \(n\times m\)
matrices which have no zero rows. Two matrices
\(G^{\prime},G^{\prime\prime}\in\mathbb{C}^{n\times m}_{\ast,r}\)
 are declared to be equivalent:
 \(G^{\prime}\stackrel{r}{\sim}
G^{\prime\prime}\), if
\begin{equation}
\label{Eqr}%
 G^{\prime\prime}=D_r G^{\prime} ,
\end{equation}
where \(D_r\) is a diagonal invertible matrix.\\
{The space} \(\mathbb{P}_{r}^{n\times (m-1)}\) is a factor-set of the set
\(\mathbb{C}^{n\times m}_{\ast,r}\) modulo the equivalency relation
\(\stackrel{r}{\sim}\).
\end{enumerate}%
\noindent The factor spaces \(\mathbb{P}_{c}^{(m-1)\times n}\) and %
\(\mathbb{P}_{r}^{n\times (m-1)}\) inherit topology from the spaces
\(\mathbb{C}^{m\times n}_{\ast,c}\) and \(\mathbb{C}^{n\times m}_{\ast,r}\),
respectively. They can be provided naturally with the structure of complex
manifolds.
\end{definition}%

If \(F^{\prime}\) and \(F^{\prime\prime}\) are two
\(\stackrel{\textup{c}}{\sim}\)\,- equivalent \(m\times n\) matrices, and
\(G^{\prime}\) and \(G^{\prime\prime}\) are two \(\stackrel{\textup{r}}{\sim}\)\,-
equivalent \(n\times m\) matrices, then the solutions \(S^{\prime}\),
\(S^{\prime\prime}\) of the Lyapunov equation \eqref{ZPLE}
 with
\(F^{\prime},G^{\prime}\) \(F^{\prime\prime},G^{\prime\prime}\), substituted instead
of  \(F,G\), and the same  \(A_{\mathcal{P}},A_{\mathcal{Z}}\) are related by
\begin{equation}%
\label{SEq} S^{\prime\prime}=D_r S^{\prime} D_c,
\end{equation}%
where \(D_c,D_r\) are the invertible diagonal matrices, which appear in \eqref{Eqc},
\eqref{Eqr}. Similar result holds also for the Lyapunov equations \eqref{PZLE}.

However, since diagonal matrices commute, the expressions on the right-hand side of
\eqref{jointrep10} will not be changed if we replace the matrices \(F,G,S\) with the
matrices \(F D_c,\, D_r G,\, D_r S D_c\), respectively.

Thus, the following result holds:

\begin{proposition}%
\label{NIBEC}%
Given \(A_{\mathcal{P}}\) and \(A_{\mathcal{Z}}\), solution of Problem \ref{pres1}
depends not on the matrices \(F,G\) themselves but on their equivalency classes in
\(\mathbb{P}_{c}^{(m-1)\times n}\), \(\mathbb{P}_{r}^{n\times (m-1)}\).
\end{proposition}

\begin{remark}%
\label{Warn}%
 In view of Remark \ref{gauge}, if \(R\) is a generic rational matrix function
  then its left and right pole semiresidual
matrix \(F_{\mathcal{P}}\) and \(G_{\mathcal{P}}\)  can be considered {\em
separately} as elements of the sets
\(\mathbb{P}_{c}^{(m-1)\times n}\) and %
\(\mathbb{P}_{r}^{n\times (m-1)}\), respectively. However, {\em simultaneously} the
matrices \(F_{\mathcal{P}}\) and \(G_{\mathcal{P}}\) can not be considered so.
The same holds for the pair of the zero semiresidual matrices, as well.
\end{remark}%
\section{Holomorphic families of generic rational matrix functions}
\label{HFGRF}
\begin{definition}\label{family}
Let \(\boldsymbol{\mathcal{D}}\) be a domain \footnote{One can also
consider a Riemann domain over \(\mathbb{C}^{2n}_*\) (see Definition
5.4.4 in \cite{Her}). }
 in the space
\(\mathbb{C}^{2n}_* \)  and for every
\(\boldsymbol{t}=(t_1,\ldots,t_{2n})\in\boldsymbol{\mathcal{D}}\) let
\(R(z,\boldsymbol{t})\) be a generic \(\mathbb{C}^{m\times m}\)-valued  rational
function of \(z\) with the pole and zero matrices
\begin{equation}\label{polematrt}
A_{\mathcal{P}}(\boldsymbol{t})=\diag(t_1,\ldots,t_n),
\quad A_{\mathcal{Z}}(\boldsymbol{t})=\diag(t_{n+1},\ldots,t_{2n}).
\end{equation}
 Assume that for every
\(\boldsymbol{t}^0\in\boldsymbol{\mathcal{D}}\) and for every fixed \(z\in
{\mathbb{C}}\setminus\{t_1^0,\ldots,t_{2n}^{0}\}\)  the function
\(R(z,\boldsymbol{t})\) is holomorphic with respect to \(\boldsymbol{t}\) in a
neighborhood of \(\boldsymbol{t}^0\). Assume also that
\begin{equation}\label{normt}
R(\infty,\boldsymbol{t})\equiv I.\end{equation} Then  the family
\(\{R(z,\boldsymbol{t})\}_{\boldsymbol{t}\in\boldsymbol{\mathcal{D}}}\) is said to
be a {\em normalized holomorphic} family of generic rational functions parameterized
by the pole and zero loci.
\end{definition}
Given  a normalized holomorphic family
\(\{R(z,\boldsymbol{t})\}_{\boldsymbol{t}\in\boldsymbol{\mathcal{D}}}\) of generic
rational functions parameterized by the pole and zero loci, we can write for each
fixed \(\boldsymbol{t}\in\boldsymbol{\mathcal{D}}\) the following representations
for  the  functions \(R(z,\boldsymbol{t})\), \(R^{-1}(z,\boldsymbol{t})\) and  the
logarithmic derivative
\begin{equation}\label{parlogder}
Q_R(z,\boldsymbol{t})\stackrel{\textup{\tiny def}}{=}\frac{\partial
R(z,\boldsymbol{t})}{\partial z}R^{-1}(z,\boldsymbol{t})\end{equation}
 (see \eqref{AELD}, \eqref{aer}):
\begin{subequations}\label{adext}
\begin{align}
R(z,\boldsymbol{t})&=I+\sum_{k=1}^n\frac{R_{k}(\boldsymbol{t})}{z-t_k},
\label{AddReprTDir}\\
R^{-1}(z,\boldsymbol{t})&=I+\sum_{k=n+1}^{2n} \frac{R_{k}(\boldsymbol{t})}{z-t_k},
\label{AddReprTInv}\\
Q_R(z,\boldsymbol{t})&= \sum_{k=1}^{2n}\frac{Q_k(\boldsymbol{t})}{z-t_k}.
\label{AELDT}
\end{align}
\end{subequations}
The residues \(R_{k}(\boldsymbol{t})\), \(Q_k(\boldsymbol{t})\),
considered as functions of \(\boldsymbol{t}\), are defined
in the whole domain \(\boldsymbol{\mathcal{D}}\). It is not hard to see
that these functions are holomorphic in \(\boldsymbol{\mathcal{D}}\):
\begin{lemma}%
\label{HRM} Let \(\boldsymbol{\mathcal{D}}\) be a domain in \(\mathbb{C}^{2n}_*\)
and let \(\{R(z,\boldsymbol{t})\}_{\boldsymbol{t}\in\boldsymbol{\mathcal{D}}}\) be a
normalized holomorphic family of generic rational functions, parameterized by the
pole and zero loci. For each fixed \(\boldsymbol{t}\in\boldsymbol{\mathcal{D}}\)
and \(1\leq k\leq n\) (respectively, \(n+1\leq k\leq 2n\))
let \(R_k(\boldsymbol{t})\) be the residue of the  rational function
 \(R(\cdot,\boldsymbol{t})\) (respectively, \(R^{-1}(\cdot,\boldsymbol{t})\))
at its pole \(t_k.\) Likewise,
for each fixed \(\boldsymbol{t}\in\boldsymbol{\mathcal{D}}\) and
\(1\leq k\leq 2n\)  let \(Q_k(\boldsymbol{t})\)  be the residue of
 the logarithmic derivative
\(Q_R(\cdot,\boldsymbol{t})\) at its pole \(t_k.\) Then  \(R_k(\boldsymbol{t})\),
\(Q_k(\boldsymbol{t})\) considered as functions of \(\boldsymbol{t}\)
are holomorphic
in \(\boldsymbol{\mathcal{D}}\).
\end{lemma}

\begin{proof}
Let us choose an arbitrary \(\boldsymbol{t}^0\in\boldsymbol{\mathcal D}\) and
\(n\) pairwise distinct points \(z_1,\ldots,z_n\) in \(\mathbb{C}
\setminus\{t_1^0,\ldots,t_{n}^{0}\}\). From the expansion
 \eqref{AddReprTDir} we derive the
following system of linear equations with respect to the residue matrices
\(R_k(\boldsymbol{t})\):
\begin{equation}%
\label{LSFRM1}%
\sum_{k=1}^{n}\frac{R_{k}(\boldsymbol{t})}{z_\ell-t_k}=R(z_\ell,\boldsymbol{t})-I,
\quad
\ell=1,\dots,n.
\end{equation}%
The matrices \(R(z_\ell,\boldsymbol{t})-I\) on the right-hand side of the system
\eqref{LSFRM1} are holomorphic with respect to \(\boldsymbol{t}\) in a neighborhood
of \(\boldsymbol{t}^0\). The determinant  of this linear system
\[\Delta(\boldsymbol{t})=\det\left(\dfrac{1}{z_\ell-t_k}\right)_{1\leq
\ell,k\leq n}\] is holomorphic in a neighborhood of \(\boldsymbol{t}^0\),
as well. In fact,
the determinant \(\Delta(\boldsymbol{t})\) (known as the
{\em Cauchy determinant}) can be calculated explicitly (see,
for example, [\cite{PS}, part VII, Section 1, No.\,3]):
\begin{equation*}%
\label{EEFD}%
\Delta(\boldsymbol{t})=\pm\frac{\prod\limits_{1\leq p<q\leq n}(z_p-z_q)(t_p-t_q)}
{\prod\limits_{\ell,k=1}^{n}(z_\ell-t_k)}.
\end{equation*}%
In particular, \(\Delta(\boldsymbol{t}^0)\not=0\). Hence, for \(k=1,\dots,n\),
  the functions
\(R_{k}(\boldsymbol{t})\)  are holomorphic
in a neighborhood of \(\boldsymbol{t}^0\).
Since this is true for any \(\boldsymbol{t}^0\in\boldsymbol{\mathcal D}\),
these functions are holomorphic in the whole domain \(\boldsymbol{\mathcal D}\).
The proof for \(R_k(\boldsymbol{t})\), when
\(n+1\leq k\leq 2n\), and for \(Q_k(\boldsymbol{t})\) is
completely analogous.
\end{proof}
\begin{remark}
\label{Hart} %
Note that, on the one hand, the
functions \(R(z,\,\boldsymbol{t})\),
\(R^{-1}(z,\boldsymbol{t})\) are rational with respect
to \(z\), and hence are holomorphic with respect to \(z\) in
\(\mathbb{C}\setminus\{t_1,\dots,t_{2n}\}\). On the other hand,
 for every fixed
\(z\in\mathbb{C}\), these  functions are holomorphic with respect to
\(\boldsymbol{t}\) in \(\boldsymbol{\mathcal
D}\setminus\bigcup_{k=1}^{2n}\{\boldsymbol{t}:t_k=z\}\). Thus, by
Hartogs theorem (see for example \textup{\cite{Shab}, Chapter I,
sections 2.3, 2.6}), the matrix functions \(R(z,\boldsymbol{t})\),
\(R^{-1}(z,\boldsymbol{t})\) are {\em jointly} holomorphic in the
variables \(z,\boldsymbol{t}\) outside the singular set
\(\bigcup_{k=1}^{2n}\{x,\boldsymbol{t}:t_k=x\}\). In view of Lemma
\ref{HRM}, the same conclusion follows
 from the representations
\eqref{AddReprTDir}, \eqref{AddReprTInv}.
\end{remark}

In order to employ the joint system representation techniques described in Sections
\ref{RMFGP} and \ref{RFPD} in the present setting, we have to establish   the
holomorphy  not only of the residues but also of the semiresidues of
\(R(\cdot,\boldsymbol{t})\). In view of Remarks \ref{gause} and \ref{gauge}, we have
a certain freedom in definition of the semiresidues and the semiresidual matrices.
 Thus one should take care in choosing  the
semiresidual matrices of \(R(\cdot,\boldsymbol{t})\) for each fixed
\(\boldsymbol{t}\) in order to obtain holomorphic functions of \(\boldsymbol{t}\).
In general,  it is possible to define the holomorphic semiresidues
 only locally
 (we refer the reader to  Appendix \ref{app2} of the present paper where the global holomorphic factorization
of a matrix function  of rank one is discussed).
\begin{lemma}%
\label{RankOneMaLe}%
Let \(M(\boldsymbol{t})\) be a \(\mathbb{C}^{m\times m}\)-valued  function,
holomorphic in a domain \({\mathcal{D}}\) of \(\mathbb{C}^N\), and let
\begin{equation}%
\label{OneRankCond}%
\rank M(\boldsymbol{t})=1\quad \forall \boldsymbol{t}\in\boldsymbol{\mathcal{D}}.
\end{equation}%
Then there exist a finite open covering \(\{\mathcal{U}_p\}_{p=1}^m\) %
of  \(\boldsymbol{\mathcal{D}}\), a collection \(\{f_p(\boldsymbol{t})\}_{p=1}^m\)
of \(\mathbb{C}^{m\times 1}\)-valued functions and a collection \(
\{g_p(\boldsymbol{t})\}_{p=1}^m\) of \(\mathbb{C}^{1\times m}\)-valued functions
satisfying the following conditions.
\begin{enumerate}
\item
For \(p=1,\dots,m\) the  functions
\(f_p(\boldsymbol{t})\) and \(g_p(\boldsymbol{t})\) are holomorphic
in \(\mathcal{U}_p\).

\item
Whenever \(\mathcal{U}_{p^\prime}\cap\mathcal{U}_{p^{\prime\prime}}\not=\emptyset\),
there exists a (scalar) function \(\varphi_{p^\prime,p^{\prime\prime}}(\boldsymbol{t})\),
holomorphic and invertible in
\(\mathcal{U}_{p^\prime}\cap\mathcal{U}_{p^{\prime\prime}}\),
such that
for every \(\boldsymbol{t}\in\mathcal{U}_{p^\prime}\cap\mathcal{U}_{p^{\prime\prime}}\)
\begin{equation}\label{gauau}
f_{p^{\prime\prime}}(\boldsymbol{t})=f_{p^{\prime}}(\boldsymbol{t})
\varphi_{p^\prime,p^{\prime\prime}}(\boldsymbol{t}),\quad
g_{p^{\prime\prime}}(\boldsymbol{t})=
\varphi^{-1}_{p^\prime,p^{\prime\prime}}(\boldsymbol{t})g_{p^{\prime}}(\boldsymbol{t}).
\end{equation}

\item For \(p=1,\dots,m\) the function
\(M(\boldsymbol{t})\) admits the factorization
\begin{equation}%
\label{LocFact}%
M(\boldsymbol{t})=f_p(\boldsymbol{t})g_p(\boldsymbol{t}),\quad
\boldsymbol{t}\in\boldsymbol{\mathcal{U}}_p.
\end{equation}%
\end{enumerate}
\end{lemma}%

\begin{proof}
Let  \(f_p(\boldsymbol{t})\) be the \(p\)-th column of the matrix
\(M(\boldsymbol{t})\) and let
\begin{equation}
\label{DefChart}
\mathcal{U}_p=\{\boldsymbol{t}\in\boldsymbol{\mathcal{D}}:f_p(\boldsymbol{t})\not=0\},\quad
p=1,\dots,m.
\end{equation}
Then, in view of \eqref{OneRankCond} and \eqref{DefChart},
\(\{\mathcal{U}_p\}_{p=1}^m\) is an open covering of  \(\boldsymbol{\mathcal{D}}\).
Furthermore, from  \eqref{OneRankCond} and \eqref{DefChart} it follows that for
\(1\leq p,q\leq m\) there exists a unique (scalar) function
\(\varphi_{p,q}(\boldsymbol{t})\), holomorphic in \(\mathcal{U}_p\), such that
\begin{equation}%
\label{connection}
f_q(\boldsymbol{t})=f_p(\boldsymbol{t})\varphi_{p,q}(\boldsymbol{t}), \quad
 \boldsymbol{t}\in \mathcal{U}_p.
\end{equation}%
Now define \(g_p(\boldsymbol{t})\) as
\begin{equation}%
g_p(\boldsymbol{t})=\begin{pmatrix}\varphi_{p,1} (\boldsymbol{t}) & \dots &
\varphi_{p,m}(\boldsymbol{t})\end{pmatrix},\quad \boldsymbol{t}\in \mathcal{U}_p.
\end{equation}%
Then the  function \(g_p(\boldsymbol{t})\) is holomorphic in \(U_p\) and, according
to  \eqref{connection}, the factorization \eqref{LocFact} holds for these
\(f_p(\boldsymbol{t})\) and \(g_p(\boldsymbol{t})\). \eqref{connection} also implies
that whenever
\(\mathcal{U}_{p^\prime}\cap\mathcal{U}_{p^{\prime\prime}}\not=\emptyset\) we have
\[
\varphi_{p^{\prime\prime},p^\prime}(\boldsymbol{t})
\varphi_{p^{\prime},k}(\boldsymbol{t})=
\varphi_{p^{\prime\prime},k}(\boldsymbol{t}),\quad
\boldsymbol{t}\in\mathcal{U}_{p^\prime}\cap\mathcal{U}_{p^{\prime\prime}},\,
1\leq k\leq m.\]
In particular,
\[\varphi_{p^{\prime\prime},p^\prime}(\boldsymbol{t})
\varphi_{p^{\prime},p^{\prime\prime}}(\boldsymbol{t})=1,\quad
\boldsymbol{t}\in\mathcal{U}_{p^\prime}\cap\mathcal{U}_{p^{\prime\prime}},\]
and \eqref{gauau} follows.
\end{proof}

\begin{theorem}\label{yokbaba1}
Let  \(\boldsymbol{\mathcal{D}}\) be a domain in \(\mathbb{C}^{2n}_*\) and let
\(\{R(z,\boldsymbol{t})\}_{\boldsymbol{t}\in\boldsymbol{\mathcal{D}}}\) be a
normalized holomorphic family of \(\mathbb{C}^{m\times m}\)-valued generic rational
functions parameterized by the pole and zero loci. Then there exist a finite open
covering\footnote{ The index \(\alpha\) runs over a finite indexing set
\(\mathfrak{A}\).} \(\{\mathcal{D}_{\alpha}\}_{\alpha\in\mathfrak{A}}\) of
 \(\boldsymbol{\mathcal{D}}\), two collections
 \(\{F_{\mathcal{P},\alpha}(\boldsymbol{t})\}_{\alpha\in\mathfrak{A}},\)
 \(\{F_{\mathcal{Z},\alpha}(\boldsymbol{t})\}_{\alpha\in\mathfrak{A}}
 \) of \(\mathbb{C}^{m\times n}\)-valued functions and two collections
 \(\{G_{\mathcal{P},\alpha}(\boldsymbol{t})\}_{\alpha\in\mathfrak{A}},\)
 \(\{G_{\mathcal{Z},\alpha}(\boldsymbol{t})\}_{\alpha\in\mathfrak{A}}
 \) of \(\mathbb{C}^{n\times m}\)-valued functions satisfying the following
 conditions.
 \begin{enumerate}
 \item For each \(\alpha\in\mathfrak{A}\) the functions
 \(F_{\mathcal{P},\alpha}(\boldsymbol{t}),\,F_{\mathcal{Z},\alpha}(\boldsymbol{t}),\,
G_{\mathcal{P},\alpha}(\boldsymbol{t}),\,G_{\mathcal{Z},\alpha}(\boldsymbol{t})\)
are holomorphic in \(\mathcal{D}_{\alpha}\). \item Whenever
\(\mathcal{D}_{\alpha^\prime}\cap\mathcal{D}_{\alpha^{\prime\prime}}\not=\emptyset\),
there exist diagonal matrix functions
\({D}_{\mathcal{P},\alpha^\prime,\alpha^{\prime\prime}}(\boldsymbol{t}),\) \(
{D}_{\mathcal{Z},\alpha^\prime,\alpha^{\prime\prime}}(\boldsymbol{t})\), holomorphic
and invertible in
\(\mathcal{D}_{\alpha^\prime}\cap\mathcal{D}_{\alpha^{\prime\prime}}\), such that
for every
\(\boldsymbol{t}\in\mathcal{D}_{\alpha^\prime}\cap\mathcal{D}_{\alpha^{\prime\prime}}\)
\begin{subequations}\label{gaugeat}
\begin{align} F_{\mathcal{P},\alpha^{\prime\prime}}(\boldsymbol{t})
= F_{\mathcal{P},\alpha^\prime}(\boldsymbol{t})
D_{\mathcal{P},\alpha^\prime,\alpha^{\prime\prime}}(\boldsymbol{t}), &\quad
G_{\mathcal{P},\alpha^{\prime\prime}}(\boldsymbol{t})=
D_{\mathcal{P},\alpha^\prime,\alpha^{\prime\prime}}^{-1}(\boldsymbol{t})
G_{\mathcal{P},\alpha^\prime}(\boldsymbol{t}), \label{polegaugea}
\\
F_{\mathcal{Z},\alpha^{\prime\prime}}(\boldsymbol{t})=
F_{\mathcal{Z},\alpha^\prime}(\boldsymbol{t})
D_{\mathcal{Z},\alpha^\prime,\alpha^{\prime\prime}}(\boldsymbol{t}), &\quad
G_{\mathcal{Z},\alpha^{\prime\prime}}(\boldsymbol{t})=
D_{\mathcal{Z},\alpha^\prime,\alpha^{\prime\prime}}^{-1}(\boldsymbol{t})
G_{\mathcal{Z},\alpha^\prime}(\boldsymbol{t}). \label{zerogaugea}
 \end{align}
 \end{subequations}
 \item For each \(\alpha\in\mathfrak{A}\) and
 \(\boldsymbol{t}\in\mathcal{D}_{\alpha}\) the matrices
\(F_{\mathcal{P},\alpha}(\boldsymbol{t})\),
\(F_{\mathcal{Z},\alpha}(\boldsymbol{t})\),
\(G_{\mathcal{P},\alpha}(\boldsymbol{t})\),
\(G_{\mathcal{Z},\alpha}(\boldsymbol{t})\)
are, respectively, the left pole, left zero, right pole, right zero
semiresidual matrices of the generic rational
function \(R(\cdot,\boldsymbol{t})\), i.e. the representations
\begin{subequations}
\label{aemfa}
\begin{align}%
\label{AddDirMaFa}%
R(z,\boldsymbol{t})&=I+F_{\mathcal{P},\alpha}(\boldsymbol{t})
\left(zI-A_{\mathcal{P}}(\boldsymbol{t})\right)^{-1}
G_{\mathcal{P},\alpha}(\boldsymbol{t}),\\
\label{AddInvMaFa}%
R^{-1}(z,\boldsymbol{t})&=I+F_{\mathcal{Z},\alpha}(\boldsymbol{t})\left(zI-A_{\mathcal{Z}}
(\boldsymbol{t})\right)^{-1}
G_{\mathcal{Z},\alpha}(\boldsymbol{t})
\end{align}%
\end{subequations}
hold true for all \(\boldsymbol{t}\in\mathcal{D}_{\alpha}\).
 \end{enumerate}
\end{theorem}

\begin{proof}
Let \(R_k(\boldsymbol{t})\), \(k=1,\dots,2n\), be the holomorphic
residue functions
as in Lemma \ref{HRM}. Since for every fixed \(\boldsymbol{t}\)
the rational function \(R(\cdot,\boldsymbol{t})\) is generic, each matrix
 \(R_k(\boldsymbol{t})\) is of rank one. Hence
there exists a finite open covering \(\{\mathcal{U}_{k,p}\}_{p=1}^m\) of
\(\boldsymbol{\mathcal{D}}\), such that in each open set \(\mathcal{U}_{k,p}\) the
function \(R_k(\boldsymbol{t})\) admits the factorization
\(R_k(\boldsymbol{t})=f_{k,p}(\boldsymbol{t})g_{k,p}(\boldsymbol{t})\) as in Lemma
\ref{RankOneMaLe}. Now it suffices to define \(\mathfrak{A}\) as the set of
\(2n\)-tuples \((p_1,\dots,p_{2n})\) such that
\(\cap_{k=1}^{2n}\mathcal{U}_{k,p_k}\not=\emptyset,\) the open covering
\(\{\mathcal{D}_{\alpha}\}_{\alpha\in\mathfrak{A}}\) of
 \(\boldsymbol{\mathcal{D}}\) by
\[\mathcal{D}_{(p_1,\dots,p_{2n})}=\bigcap_{k=1}^{2n}\mathcal{U}_{k,p_k},\]
and the collections
 \(\{F_{\mathcal{P},\alpha}(\boldsymbol{t})\}_{\alpha\in\mathfrak{A}},\)
\(\{F_{\mathcal{Z},\alpha}(\boldsymbol{t})\}_{\alpha\in\mathfrak{A}}
 \),
 \(\{G_{\mathcal{P},\alpha}(\boldsymbol{t})\}_{\alpha\in\mathfrak{A}},\)
 \(\{G_{\mathcal{Z},\alpha}(\boldsymbol{t})\}_{\alpha\in\mathfrak{A}}
 \) by
\begin{align*}
F_{\mathcal{P},(p_1,\dots,p_{2n})}(\boldsymbol{t})& = \begin{pmatrix}f_{1,p_1}(\boldsymbol{t})
& \dots & f_{n,p_n}(\boldsymbol{t})\end{pmatrix} ,\\
F_{\mathcal{Z},(p_1,\dots,p_{2n})}(\boldsymbol{t})
& = \begin{pmatrix}f_{n+ 1,p_{n+1}}(\boldsymbol{t}) & \dots & f_{2n,p_{2n}}(\boldsymbol{t})
\end{pmatrix} ,\\[1ex]
G_{\mathcal{P},(p_1,\dots,p_{2n})}(\boldsymbol{t})& =\begin{pmatrix} g_{1,p_1}(\boldsymbol{t})
  \\ \vdots \\ g_ {n,p_n}(\boldsymbol{t})\end{pmatrix} ,\quad
G_{\mathcal{Z},(p_1,\dots,p_{2n})}(\boldsymbol{t}) =
\begin{pmatrix}g_ {n+ 1,p_{n+1}}(\boldsymbol{t})\\ \vdots \\
 g_{2n,p_{2n}}(\boldsymbol{t})\end{pmatrix}.
\end{align*}
\end{proof}

\begin{definition}
Let  \(\boldsymbol{\mathcal{D}}\) be a domain in \(\mathbb{C}^{2n}_*\) and let
\(\{R(z,\boldsymbol{t})\}_{\boldsymbol{t}\in\boldsymbol{\mathcal{D}}}\) be a
normalized holomorphic family of \(\mathbb{C}^{m\times m}\)-valued generic rational
functions parameterized by the pole and zero loci. Let a finite open covering
\(\{\mathcal{D}_{\alpha}\}_{\alpha\in\mathfrak{A}}\) of
 \(\boldsymbol{\mathcal{D}}\),  collections
 \(\{F_{\mathcal{P},\alpha}(\boldsymbol{t})\}_{\alpha\in\mathfrak{A}},\)
 \(\{F_{\mathcal{Z},\alpha}(\boldsymbol{t})\}_{\alpha\in\mathfrak{A}}
 \) of \(\mathbb{C}^{m\times n}\)-valued functions and collections
 \(\{G_{\mathcal{P},\alpha}(\boldsymbol{t})\}_{\alpha\in\mathfrak{A}},\)
 \(\{G_{\mathcal{Z},\alpha}(\boldsymbol{t})\}_{\alpha\in\mathfrak{A}}
 \) of \(\mathbb{C}^{n\times m}\)-valued functions
satisfy the conditions 1. -- 3.  of Theorem \ref{yokbaba1}.
Then the collections \(\{F_{\mathcal{P},\alpha}(\boldsymbol{t})\}_{\alpha\in\mathfrak{A}},\)
 \(\{F_{\mathcal{Z},\alpha}(\boldsymbol{t})\}_{\alpha\in\mathfrak{A}}
 \),
\(\{G_{\mathcal{P},\alpha}(\boldsymbol{t})\}_{\alpha\in\mathfrak{A}},\)
 \(\{G_{\mathcal{Z},\alpha}(\boldsymbol{t})\}_{\alpha\in\mathfrak{A}}
 \) are said to be the collections of, respectively,
the {\em left pole, left zero, right pole, right zero semiresidual
functions} related to the family
\(\{R(z,\boldsymbol{t})\}_{\boldsymbol{t}\in\boldsymbol{\mathcal{D}}}\).
\end{definition}

Now
 we can tackle the
problem of recovery of a holomorphic family of generic matrix functions from
the semiresidual data.
Once again, let \(\boldsymbol{\mathcal{D}}\) be a domain in \(\mathbb{C}^{2n}_*\) and let
\(\{R(z,\boldsymbol{t})\}_{\boldsymbol{t}\in\boldsymbol{\mathcal{D}}}\) be a
normalized holomorphic family of \(\mathbb{C}^{m\times m}\)-valued generic rational
functions, parameterized by the pole and zero loci. Let  collections
 \(\{F_{\mathcal{P},\alpha}(\boldsymbol{t})\}_{\alpha\in\mathfrak{A}},\)
 \(\{F_{\mathcal{Z},\alpha}(\boldsymbol{t})\}_{\alpha\in\mathfrak{A}}
 \),
 \(\{G_{\mathcal{P},\alpha}(\boldsymbol{t})\}_{\alpha\in\mathfrak{A}},\)
 \(\{G_{\mathcal{Z},\alpha}(\boldsymbol{t})\}_{\alpha\in\mathfrak{A}}
 \) be the collections of the semiresidual matrices related to the family
\(\{R(z,\boldsymbol{t})\}_{\boldsymbol{t}\in\boldsymbol{\mathcal{D}}}\).
Then for each \(\alpha\in\mathfrak{A}\) and a fixed \(\boldsymbol{t}\in\mathcal{D}_{\alpha}\)
the matrices \(A_{\mathcal{P}}(\boldsymbol{t})\),
 \(A_{\mathcal{Z}}(\boldsymbol{t})\)
are the pole and zero matrices of the generic rational function
\(R(z,\boldsymbol{t})\), and \(F_{\mathcal{P},\alpha}(\boldsymbol{t})\),
\(F_{\mathcal{Z},\alpha}(\boldsymbol{t})\),
\(G_{\mathcal{P},\alpha}(\boldsymbol{t})\),
\(G_{\mathcal{Z},\alpha}(\boldsymbol{t})\) are, respectively, the left pole, left
zero, right pole, right zero semiresidual matrices of \(R(\cdot,\boldsymbol{t})\).
The appropriate coupling matrices \( S_{\mathcal{ZP},\alpha}(\boldsymbol{t})\), \(
S_{\mathcal{PZ},\alpha}(\boldsymbol{t})\) satisfying the Lyapunov equations
\begin{subequations}\label{SLt}
\begin{align}
A_{\mathcal{Z}}(\boldsymbol{t})S_{\mathcal{ZP},\alpha}(\boldsymbol{t})-S_{\mathcal{ZP},\alpha}
(\boldsymbol{t})
A_{\mathcal{P}}(\boldsymbol{t})&=G_{\mathcal{Z},\alpha}(\boldsymbol{t})F_{\mathcal{P},\alpha}
(\boldsymbol{t}), \label{SL3t}\\
A_{\mathcal{P}}(\boldsymbol{t})S_{\mathcal{PZ},\alpha}(\boldsymbol{t})-S_{\mathcal{PZ},\alpha}
(\boldsymbol{t}) A_{\mathcal{Z}}(\boldsymbol{t})&=G_{\mathcal
{P},\alpha}(\boldsymbol{t})F_{\mathcal{Z},\alpha}(\boldsymbol{t})
\label{SL4t}%
\end{align}
\end{subequations}
are  given  by
\begin{equation}%
\label{ESNPtt}%
S_{\mathcal{ZP},\alpha}(\boldsymbol{t})
=\left(\frac{g_{{n+i},\alpha}(\boldsymbol{t})f_{j,\alpha}(\boldsymbol{t})}
{t_{n+i}-t_j} \right)_{i,j=1}^n,\quad
S_{\mathcal{PZ},\alpha}(\boldsymbol{t})=
\left(\frac{g_{i,\alpha}(\boldsymbol{t})f_{n+j,\alpha}(\boldsymbol{t})}{t_i-t_{n+j}}
\right)_{i,j=1}^n,
\end{equation}%
where for \(1\leq k\leq n\) \(g_{k,\alpha}(\boldsymbol{t})\) is the \(k\)-th row of
\(G_{\mathcal {P},\alpha}(\boldsymbol{t})\), \(g_{n+k,\alpha}(\boldsymbol{t})\))
 is the \(k\)-th row of \(G_{\mathcal {Z},\alpha}(\boldsymbol{t})\),
\(f_{k,\alpha}(\boldsymbol{t})\) is the \(k\)-th column of \(F_{\mathcal
{P},\alpha}(\boldsymbol{t})\),  \(f_{n+k,\alpha}(\boldsymbol{t})\)) is the \(k\)-th
column of  \(F_{\mathcal {Z},\alpha}(\boldsymbol{t})\) (compare with similar
expressions  \eqref{leftsemires}, \eqref{rightsemires}, \eqref{ESNP}). From the
explicit expressions \eqref{ESNPtt} it is evident that \(
S_{\mathcal{ZP},\alpha}(\boldsymbol{t})\), \(
S_{\mathcal{PZ},\alpha}(\boldsymbol{t})\), considered as functions of
\(\boldsymbol{t}\), are holomorphic in \(\mathcal{D}_{\alpha}\). According to
Proposition \ref{MuInc}, the functions \( S_{\mathcal{ZP},\alpha}(\boldsymbol{t})\),
\( S_{\mathcal{PZ},\alpha}(\boldsymbol{t})\) are mutually inverse:
\begin{equation}%
\label{MIRtt}%
S_{\mathcal{ZP},\alpha}(\boldsymbol{t}) S_{\mathcal{PZ},\alpha}(\boldsymbol{t})
=S_{\mathcal{PZ},\alpha}(\boldsymbol{t}) S_{\mathcal{ZP},\alpha}(\boldsymbol{t})=I,
\quad \boldsymbol{t}\in\mathcal{D}_{\alpha},
\end{equation}
and
 the following relations hold:
\begin{equation}
\label{ZPCRtt}%
G_{\mathcal{Z},\alpha}(\boldsymbol{t})=-S_{\mathcal{ZP},\alpha}(\boldsymbol{t})
G_{\mathcal{P},\alpha}(\boldsymbol{t}),\quad F_{\mathcal{P},\alpha}(\boldsymbol{t})=
F_{\mathcal{Z},\alpha}(\boldsymbol{t})S_{\mathcal{ZP},\alpha}(\boldsymbol{t}).
\end{equation}
Furthermore, for  \(\boldsymbol{t}\in\mathcal{D}_{\alpha}\) the function
\(R(z,\boldsymbol{t})\) admits the  representation
\begin{multline}\label{jointrep0tt}
R(z,\boldsymbol{t})R^{-1}(\omega,\boldsymbol{t})
=I+\\+(z-\omega)F_{\mathcal{P},\alpha}(\boldsymbol{t})
(zI-A_{\mathcal{P}}(\boldsymbol{t}))^{-1}S_{\mathcal{ZP},\alpha}^{-1}(\boldsymbol{t})
(\omega
I-A_{\mathcal{Z}}(\boldsymbol{t}))^{-1}G_{\mathcal{Z},\alpha}(\boldsymbol{t}),\end{multline}
and, in view of Proposition \ref{ReprLogDer}, its logarithmic derivative with
respect to \(z\) admits the representation
\begin{equation}
\dfrac{\partial R(z,\boldsymbol{t})}{\partial z}R^{-1}(z,\boldsymbol{t})
= F_{\mathcal{P},\alpha}(\boldsymbol{t})
 (zI-A_{\mathcal{P}}(\boldsymbol{t}))^{-1}S_{\mathcal{ZP},\alpha}^{-1}(\boldsymbol{t})
(zI-A_{\mathcal{Z}}(\boldsymbol{t}))^{-1}G_{\mathcal{Z},\alpha}(\boldsymbol{t}).
\label{replogdertt}
\end{equation}

The representations \eqref{jointrep0tt}, \eqref{replogdertt} above are local: each
of them holds in the appropriate individual subset \(\mathcal{D}_\alpha\). Note,
however, that whenever
\(\mathcal{D}_{\alpha^\prime}\cap\mathcal{D}_{\alpha^{\prime\prime}}\not=\emptyset\),
by Theorem \ref{yokbaba1} we have
\[S_{\mathcal{ZP},\alpha^{\prime\prime}}(\boldsymbol{t})=
D_{\mathcal{Z},\alpha^\prime,\alpha^{\prime\prime}}^{-1}(\boldsymbol{t})
S_{\mathcal{ZP},\alpha^{\prime}}(\boldsymbol{t})
D_{\mathcal{P},\alpha^\prime,\alpha^{\prime\prime}}(\boldsymbol{t}),\quad
\boldsymbol{t}\in
\mathcal{D}_{\alpha^\prime}\cap\mathcal{D}_{\alpha^{\prime\prime}},\] where  the
functions \(D_{\mathcal{Z},\alpha^\prime,\alpha^{\prime\prime}}(\boldsymbol{t})\),
\(D_{\mathcal{P},\alpha^\prime,\alpha^{\prime\prime}}(\boldsymbol{t})\) are as in
\eqref{gaugeat}. Hence the expressions \eqref{jointrep0tt}, \eqref{replogdertt}
coincide in the intersections of the subsets \(\mathcal{D}_\alpha\) (although the
individual functions \(F_{\mathcal{P},\alpha}(\boldsymbol{t})\),
\(S_{\mathcal{ZP},\alpha}(\boldsymbol{t})\),
\(G_{\mathcal{Z},\alpha}(\boldsymbol{t})\) do not). Here it is a self-evident fact,
because  these  expressions represent  globally defined objects. In the next
section, where we shall use such local representations to {\em construct} globally
defined objects, it will become a requirement.

\section{Holomorphic families of generic rational matrix
functions with prescribed local data}
\label{HFPD}

This section can be considered as a \(\boldsymbol{t}\)\,-\,dependent version of
Section \ref{RFPD}. Here we consider the problem, how to construct a normalized
holomorphic family of generic rational functions\footnote{ See Definition
\ref{family}.} with prescribed local data. The nature of such data is suggested by
the considerations of the previous section (see, in particular, Theorem
\ref{yokbaba1}).

Let \(\boldsymbol{\mathcal{D}}\) be a domain in \(\mathbb{C}^{2n}_*\) and let
   \(\{\mathcal{D}_{\alpha}\}_{\alpha\in\mathfrak{A}}\) be a finite open
covering of
 \(\boldsymbol{\mathcal{D}}\).  We assume that
 \(\{F_{\alpha}(\boldsymbol{t})\}_{\alpha\in\mathfrak{A}}\) and
\(\{G_{\alpha}(\boldsymbol{t})\}_{\alpha\in\mathfrak{A}}\) are collections
  of, respectively, \(\mathbb{C}^{m\times n}\)-valued and \(\mathbb{C}^{n\times m}\)-valued
functions, satisfying the following conditions:
\begin{enumerate}
\item For each \(\alpha\in\mathfrak{A}\) the functions
 \(F_{\alpha}(\boldsymbol{t})\),
\(G_{\alpha}(\boldsymbol{t})\)
are holomorphic in \(\mathcal{D}_{\alpha}\).

\item Whenever
\(\mathcal{D}_{\alpha^\prime}\cap\mathcal{D}_{\alpha^{\prime\prime}}\not=\emptyset\),
there exist diagonal matrix functions
\({D}_{r,\alpha^\prime,\alpha^{\prime\prime}}(\boldsymbol{t}),\) \(
{D}_{c,\alpha^\prime,\alpha^{\prime\prime}}(\boldsymbol{t})\), holomorphic
and invertible in
\(\mathcal{D}_{\alpha^\prime}\cap\mathcal{D}_{\alpha^{\prime\prime}}\), such that
for every
\(\boldsymbol{t}\in\mathcal{D}_{\alpha^\prime}\cap\mathcal{D}_{\alpha^{\prime\prime}}\)
\begin{equation}\label{pgauz}
F_{\alpha^{\prime\prime}}(\boldsymbol{t})=
F_{\alpha^\prime}(\boldsymbol{t})
D_{c,\alpha^\prime,\alpha^{\prime\prime}}(\boldsymbol{t}), \quad
G_{\alpha^{\prime\prime}}(\boldsymbol{t})=
D_{r,\alpha^\prime,\alpha^{\prime\prime}}(\boldsymbol{t})
G_{\alpha^\prime}(\boldsymbol{t}). \end{equation}
\end{enumerate}

\begin{remark}
Conditions 1. and 2. imply that the collections
\(\{F_{\alpha}(\boldsymbol{t})\}_{\alpha\in\mathfrak{A}}\) and
\(\{G_{\alpha}(\boldsymbol{t})\}_{\alpha\in\mathfrak{A}}\)
represent  holomorphic mappings
from \(\boldsymbol{\mathcal{D}}\) into
the  spaces\footnote{See Definition \ref{ProjEquivDef}.}
 \(\mathbb{P}_{c}^{(m-1)\times n}\) and
\(\mathbb{P}_{r}^{n\times (m-1)}\), respectively.
\end{remark}

For each \(\alpha\in\mathfrak{A}\) and
\(\boldsymbol{t}\in\mathcal{D}_{\alpha}\)
we consider the Lyapunov equation
\begin{equation}
\label{SL1TPr}
A_{\mathcal{P}}(\boldsymbol{t})S_{\alpha}(\boldsymbol{t})-
S_{\alpha}(\boldsymbol{t})A_{\mathcal{Z}}(\boldsymbol{t})=
G_{\alpha}(\boldsymbol{t})F_{\alpha}(\boldsymbol{t}),
\end{equation}
where \(A_{\mathcal{P}}(\boldsymbol{t})\) and \(A_{\mathcal{Z}}(\boldsymbol{t})\) are as in \eqref{polematrt}.
Its solution \(S_{\alpha}(\boldsymbol{t})\)
is the \(\mathcal{PZ}\)-coupling\footnote{See Definition \ref{DefLocData}.}
 matrix, related to the data
\( A_{\mathcal{P}}(\boldsymbol{t})\), \( A_{\mathcal{Z}}(\boldsymbol{t})\),
\(F_{\alpha}(\boldsymbol{t})\), \(G_{\alpha}(\boldsymbol{t})\).
Considered as a function of \(\boldsymbol{t}\),
\(S_{\alpha}(\boldsymbol{t})\)
is  holomorphic in
\(\mathcal{D}_{\alpha}\), because the right-hand side
\(G_{\alpha}(\boldsymbol{t})F_{\alpha}(\boldsymbol{t})\)
is  holomorphic in
\(\mathcal{D}_{\alpha}\). The collection of functions \(
\{S_{\alpha}(\boldsymbol{t})\}_{\alpha\in\mathfrak{A}}
\) is said to be the collection of \(\mathcal{PZ}\)-{\em coupling
functions} related to the pair of collections
 \(\{F_{\alpha}(\boldsymbol{t})\}_{\alpha\in\mathfrak{A}}\),
\(\{G_{\alpha}(\boldsymbol{t})\}_{\alpha\in\mathfrak{A}}\).
 In view of \eqref{pgauz}, whenever
\(\mathcal{D}_{\alpha^\prime}\cap\mathcal{D}_{\alpha^{\prime\prime}}\not=\emptyset\),
we have
\begin{equation}\label{glue1}
S_{\alpha^{\prime\prime}}(\boldsymbol{t})=
D_{r,\alpha^\prime,\alpha^{\prime\prime}}(\boldsymbol{t})
S_{\alpha^{\prime\prime}}(\boldsymbol{t})D_{c,\alpha^\prime,\alpha^{\prime\prime}}(\boldsymbol{t}),
\quad
\boldsymbol{t}\in \mathcal{D}_{\alpha^\prime}\cap\mathcal{D}_{\alpha^{\prime\prime}},
\end{equation}
where \(D_{r,\alpha^\prime,\alpha^{\prime\prime}}(\boldsymbol{t})\),
\(D_{c,\alpha^\prime,\alpha^{\prime\prime}}(\boldsymbol{t})\) are diagonal,
holomorphic and invertible matrix functions.
Hence either \(\forall \alpha\in\mathfrak{A}\) \(\det S_{\alpha}(\boldsymbol{t})\equiv 0\)
or \(\forall \alpha\in\mathfrak{A}\) \(\det S_{\alpha}(\boldsymbol{t})\not\equiv 0.\)
In the latter case the pair of
collections \(\{F_{\alpha}(\boldsymbol{t})\}_{\alpha\in\mathfrak{A}}\) and
\(\{G_{\alpha}(\boldsymbol{t})\}_{\alpha\in\mathfrak{A}}\)
is said to be \(\mathcal{PZ}\)-{\em admissible},
and the set
\begin{equation}\label{defgpz}
\Gamma_{\mathcal{PZ}}=\bigcup_{\alpha\in\mathfrak{A}}\{\boldsymbol{t}\in
\mathcal{D}_{\alpha}:
\det S_{\alpha}(\boldsymbol{t})= 0\}
\end{equation}
is said to be the \(\mathcal{PZ}\)-{\em singular} set related to
 the pair of collections \(\{F_{\alpha}(\boldsymbol{t})\}_{\alpha\in\mathfrak{A}}\),
\(\{G_{\alpha}(\boldsymbol{t})\}_{\alpha\in\mathfrak{A}}\) .

In the same way, we can consider
the collection of functions \(
\{S_{\alpha}(\boldsymbol{t})\}_{\alpha\in\mathfrak{A}}
\), where for each \(\alpha\in\mathfrak{A}\) and
\(\boldsymbol{t}\in\mathcal{D}_{\alpha}\)
 the matrix \(S_{\alpha}(\boldsymbol{t})\) is the solution
of the Lyapunov equation
\begin{equation}
\label{SL2TPr}
A_{\mathcal{Z}}(\boldsymbol{t})S_{\alpha}(\boldsymbol{t})-
S_{\alpha}(\boldsymbol{t})A_{\mathcal{P}}(\boldsymbol{t})=
G_{\alpha}(\boldsymbol{t})F_{\alpha}(\boldsymbol{t}).
\end{equation}
This collection is said to be the collection of \(\mathcal{ZP}\)-{\em coupling
functions} related to the pair of collections
 \(\{F_{\alpha}(\boldsymbol{t})\}_{\alpha\in\mathfrak{A}}\),
\(\{G_{\alpha}(\boldsymbol{t})\}_{\alpha\in\mathfrak{A}}\).
If for every \(\alpha\in\mathfrak{A}\)
 \(\det S_{\alpha}(\boldsymbol{t})\not\equiv 0\) in \(\mathcal{D}_{\alpha}\) then
the pair of collections \(\{F_{\alpha}(\boldsymbol{t})\}_{\alpha\in\mathfrak{A}}\),
\(\{G_{\alpha}(\boldsymbol{t})\}_{\alpha\in\mathfrak{A}}\) is said to be
\(\mathcal{ZP}\)-{\em admissible}, and the set
\begin{equation}\label{defgzp}
\Gamma_{\mathcal{ZP}}=\bigcup_{\alpha\in\mathfrak{A}}\{\boldsymbol{t}
\in\mathcal{D}_{\alpha}:
\det S_{\alpha}(\boldsymbol{t})= 0\}
\end{equation}
is said to be the \(\mathcal {ZP}\)-{\em singular} set related to
 the pair of collections \(\{F_{\alpha}(\boldsymbol{t})\}_{\alpha\in\mathfrak{A}}\),
\(\{G_{\alpha}(\boldsymbol{t})\}_{\alpha\in\mathfrak{A}}\).

\begin{remark}\label{sinman}
Note that, in view of \eqref{pgauz},
 if for each \(\alpha\in\mathfrak{A}\) \(S_{\alpha}(\boldsymbol{t})\) satisfies
the Lyapunov equation \eqref{SL1TPr} (or \eqref{SL2TPr}) then the subsets
\[\Gamma_\alpha=\{\boldsymbol{t}
\in\mathcal{D}_{\alpha}:
\det S_{\alpha}(\boldsymbol{t})= 0\}\]
of the appropriate singular set agree in the
intersections of the sets \(\mathcal{D}_\alpha\):
\[\Gamma_{\alpha^{\prime}}
\cap(\mathcal{D}_{\alpha^\prime}\cap\mathcal{D}_{\alpha^{\prime\prime}})=
\Gamma_{\alpha^{\prime\prime}}
\cap(\mathcal{D}_{\alpha^\prime}\cap\mathcal{D}_{\alpha^{\prime\prime}})\quad
\forall \alpha^\prime,\alpha^{\prime\prime}.
\]
\end{remark}

\begin{theorem}\label{main1}
Let \(\boldsymbol{\mathcal{D}}\) be a domain in \(\mathbb{C}^{2n}_*\) and let
   \(\{\mathcal{D}_{\alpha}\}_{\alpha\in\mathfrak{A}}\) be a finite open
covering of
 \(\boldsymbol{\mathcal{D}}\).  Let
 \(\{F_{\alpha}(\boldsymbol{t})\}_{\alpha\in\mathfrak{A}}\) and
\(\{G_{\alpha}(\boldsymbol{t})\}_{\alpha\in\mathfrak{A}}\) are collections
  of, respectively,
\(\mathbb{C}^{m\times n}\)-valued and \(\mathbb{C}^{n\times m}\)-valued
functions, satisfying the following conditions:
\begin{enumerate}
\item For each \(\alpha\in\mathfrak{A}\) the functions
 \(F_{\alpha}(\boldsymbol{t})\),
\(G_{\alpha}(\boldsymbol{t})\)
are holomorphic in \(\mathcal{D}_{\alpha}\).

\item Whenever
\(\mathcal{D}_{\alpha^\prime}\cap\mathcal{D}_{\alpha^{\prime\prime}}\not=\emptyset\),
there exist diagonal matrix functions
\({D}_{r,\alpha^\prime,\alpha^{\prime\prime}}(\boldsymbol{t}),\) \(
{D}_{c,\alpha^\prime,\alpha^{\prime\prime}}(\boldsymbol{t})\), holomorphic
and invertible in
\(\mathcal{D}_{\alpha^\prime}\cap\mathcal{D}_{\alpha^{\prime\prime}}\), such that
for every
\(\boldsymbol{t}\in\mathcal{D}_{\alpha^\prime}\cap\mathcal{D}_{\alpha^{\prime\prime}}\)
\begin{equation*}
F_{\alpha^{\prime\prime}}(\boldsymbol{t})= F_{\alpha^\prime}(\boldsymbol{t})
D_{c,\alpha^\prime,\alpha^{\prime\prime}}(\boldsymbol{t}), \quad
G_{\alpha^{\prime\prime}}(\boldsymbol{t})=
D_{r,\alpha^\prime,\alpha^{\prime\prime}}(\boldsymbol{t})
G_{\alpha^\prime}(\boldsymbol{t}). \end{equation*}

 \item  The pair of collections
\(\{F_{\alpha}(\boldsymbol{t})\}_{\alpha\in\mathfrak{A}}\),
\(\{G_{\alpha}(\boldsymbol{t})\}_{\alpha\in\mathfrak{A}}\) is  \(\mathcal{ZP}\)-{\em
admissible}. \end{enumerate}
 Let \(\Gamma_{\mathcal{ZP}}\) denote the
 \(\mathcal{ZP}\)-singular set related to the pair of collections
\(\{F_{\alpha}(\boldsymbol{t})\}_{\alpha\in\mathfrak{A}}\),
\(\{G_{\alpha}(\boldsymbol{t})\}_{\alpha\in\mathfrak{A}}\).
 Then there exists a unique normalized
holomorphic family \(\{R(z,\boldsymbol{t})\}_{\boldsymbol{t}
\in\boldsymbol{\mathcal{D}}\setminus \Gamma_{\mathcal{ZP}}}\)
 of rational generic
functions parameterized by the pole and zero loci  such that for every
\(\alpha\in\mathfrak{A}\) and
 \(\boldsymbol{t}\in\mathcal{D}_{\alpha}\setminus\Gamma_{\mathcal{ZP}}\)
the matrices \(F_{\alpha}(\boldsymbol{t})\) and \(G_{\alpha}(\boldsymbol{t})\)
are, respectively, the left pole and right zero semiresidual matrices
of \(R(\cdot,\boldsymbol{t})\):
\begin{equation}\label{kokozp}
F_{\alpha}(\boldsymbol{t})=F_{\mathcal{P},\alpha}(\boldsymbol{t}),
\quad G_{\alpha}(\boldsymbol{t})=G_{\mathcal{Z},\alpha}(\boldsymbol{t}),\quad
\forall\boldsymbol{t}\in\mathcal{D}_{\alpha}\setminus\Gamma_{\mathcal{ZP}}.
\end{equation}
It is locally given by
\begin{equation}\label{jointrep1t}
R(z,\boldsymbol{t}) =I-F_{\alpha}(\boldsymbol{t})
(zI-A_{\mathcal{P}}(\boldsymbol{t}))^{-1}S_\alpha^{-1}(\boldsymbol{t})
G_{\alpha}(\boldsymbol{t}),
\quad\boldsymbol{t}\in\mathcal{D}_{\alpha}\setminus\Gamma_{\mathcal{ZP}},\end{equation}
where \(\{S_\alpha(\boldsymbol{t})\}_{\alpha\in\mathfrak{A}}\) is the collection of
\(\mathcal{ZP}\)-coupling functions related to the pair of collections
\(\{F_{\alpha}(\boldsymbol{t})\}_{\alpha\in\mathfrak{A}}\),
\(\{G_{\alpha}(\boldsymbol{t})\}_{\alpha\in\mathfrak{A}}\). Furthermore,  the
logarithmic derivative of \(R(z,\boldsymbol{t})\) with respect to \(z\) admits the
local representation
\begin{multline}
\dfrac{\partial R(z,\boldsymbol{t})}{\partial z}R^{-1}(z,\boldsymbol{t})=\\ =
F_{\alpha}(\boldsymbol{t})
 (zI-A_{\mathcal{P}}(\boldsymbol{t}))^{-1}S_\alpha^{-1}(\boldsymbol{t})
(zI-A_{\mathcal{Z}}(\boldsymbol{t}))^{-1}G_{\alpha}(\boldsymbol{t}),\\
\boldsymbol{t}\in\mathcal{D}_{\alpha}\setminus\Gamma_{\mathcal{ZP}}.
\label{replogder1t}
\end{multline}
\end{theorem}

\begin{proof}
In view of Proposition \ref{record} and condition 3, for each
\(\alpha\in\mathfrak{A}\) and \(\boldsymbol{t}\in\mathcal{D}_{\alpha}
\setminus\Gamma_{\mathcal{ZP}}\) there exists a unique generic rational function
\(R_\alpha(\cdot,\boldsymbol{t})\), normalized by
\(R_\alpha(\infty,\boldsymbol{t})=I\), with the pole and zero matrices
\(A_{\mathcal{P}}(\boldsymbol{t})\), \(A_{\mathcal{Z}}(\boldsymbol{t})\) and the
prescribed left zero and right pole semiresidual matrices \eqref{kokozp}. The
function \(R_\alpha(\cdot,\boldsymbol{t})\) and its logarithmic derivative admit the
representations
\begin{equation}\label{jointrep1ta}
R_\alpha(z,\boldsymbol{t}) =I-F_{\alpha}(\boldsymbol{t})
(zI-A_{\mathcal{P}}(\boldsymbol{t}))^{-1}S_\alpha^{-1}(\boldsymbol{t})
G_{\alpha}(\boldsymbol{t}),\end{equation}
\begin{multline*}
R_\alpha^\prime(z,\boldsymbol{t})R_\alpha^{-1}(z,\boldsymbol{t})=\\ =
F_{\alpha}(\boldsymbol{t})
 (zI-A_{\mathcal{P}}(\boldsymbol{t}))^{-1}S_\alpha^{-1}(\boldsymbol{t})
(zI-A_{\mathcal{Z}}(\boldsymbol{t}))^{-1}G_{\alpha}(\boldsymbol{t}).
\end{multline*}
 From the representation \eqref{jointrep1ta} and condition 1
it follows that the family
\(\{R_\alpha(z,\boldsymbol{t})\}_{\boldsymbol{t}\in\boldsymbol{\mathcal{D}}_{\alpha}\setminus
\Gamma_{\mathcal{ZP}}}\) is  holomorphic.
 In view of condition 2 (see also \eqref{glue1}),
 we have
 \[R_{\alpha^{\prime}}(z,\boldsymbol{t})=R_{\alpha^{\prime\prime}}(z,\boldsymbol{t}),
 \quad \forall\boldsymbol{t}\in
 (\mathcal{D}_{\alpha^\prime}\cap\mathcal{D}_{\alpha^{\prime\prime}})
 \setminus\Gamma_{\mathcal{ZP}}.\]
Hence we can define the holomorphic family
\(\{R(z,\boldsymbol{t})\}_{\boldsymbol{t}\in\boldsymbol{\mathcal{D}}\setminus
\Gamma_{\mathcal{ZP}}}\) by
\[R(z,\boldsymbol{t})=R_\alpha(z,\boldsymbol{t}),\quad
\boldsymbol{t}\in\boldsymbol{\mathcal{D}}_{\alpha}\setminus \Gamma_{\mathcal{ZP}}\]
to obtain the local representations \eqref{jointrep1t}, \eqref{replogder1t}. The
uniqueness of such a family follows from the uniqueness of each function
\(R_\alpha(\cdot,\boldsymbol{t})\).
\end{proof}

\begin{theorem}\label{main2}
Let \(\boldsymbol{\mathcal{D}}\) be a domain in \(\mathbb{C}^{2n}_*\) and let
   \(\{\mathcal{D}_{\alpha}\}_{\alpha\in\mathfrak{A}}\) be a finite open
covering of
 \(\boldsymbol{\mathcal{D}}\).  Let
 \(\{F_{\alpha}(\boldsymbol{t})\}_{\alpha\in\mathfrak{A}}\) and
\(\{G_{\alpha}(\boldsymbol{t})\}_{\alpha\in\mathfrak{A}}\) are collections
  of, respectively,
\(\mathbb{C}^{m\times n}\)-valued and \(\mathbb{C}^{n\times m}\)-valued functions,
satisfying the following conditions:
\begin{enumerate}
\item For each \(\alpha\in\mathfrak{A}\) the functions
 \(F_{\alpha}(\boldsymbol{t})\),
\(G_{\alpha}(\boldsymbol{t})\) are holomorphic in \(\mathcal{D}_{\alpha}\).

\item Whenever
\(\mathcal{D}_{\alpha^\prime}\cap\mathcal{D}_{\alpha^{\prime\prime}}\not=\emptyset\),
there exist diagonal matrix functions
\({D}_{r,\alpha^\prime,\alpha^{\prime\prime}}(\boldsymbol{t}),\) \(
{D}_{c,\alpha^\prime,\alpha^{\prime\prime}}(\boldsymbol{t})\), holomorphic and
invertible in
\(\mathcal{D}_{\alpha^\prime}\cap\mathcal{D}_{\alpha^{\prime\prime}}\), such that
for every
\(\boldsymbol{t}\in\mathcal{D}_{\alpha^\prime}\cap\mathcal{D}_{\alpha^{\prime\prime}}\)
\begin{equation*}
F_{\alpha^{\prime\prime}}(\boldsymbol{t})= F_{\alpha^\prime}(\boldsymbol{t})
D_{c,\alpha^\prime,\alpha^{\prime\prime}}(\boldsymbol{t}), \quad
G_{\alpha^{\prime\prime}}(\boldsymbol{t})=
D_{r,\alpha^\prime,\alpha^{\prime\prime}}(\boldsymbol{t})
G_{\alpha^\prime}(\boldsymbol{t}). \end{equation*}

 \item The pair of collections
\(\{F_{\alpha}(\boldsymbol{t})\}_{\alpha\in\mathfrak{A}}\),
\(\{G_{\alpha}(\boldsymbol{t})\}_{\alpha\in\mathfrak{A}}\) is  \(\mathcal{PZ}\)-{\em
admissible}.
\end{enumerate}
Let \(\Gamma_{\mathcal{PZ}}\) denote the \(\mathcal{PZ}\)-singular set related to
the pair of collections \(\{F_{\alpha}(\boldsymbol{t})\}_{\alpha\in\mathfrak{A}}\),
\(\{G_{\alpha}(\boldsymbol{t})\}_{\alpha\in\mathfrak{A}}\).
 Then there
exists a unique normalized holomorphic family
\(\{R(z,\boldsymbol{t})\}_{\boldsymbol{t}\in\boldsymbol{\mathcal{D}}\setminus
\Gamma_{\mathcal{PZ}}}\)
 of rational generic
functions parameterized by the pole and zero loci such that for every
\(\alpha\in\mathfrak{A}\) and   \(\boldsymbol{t}\in\mathcal{D}_{\alpha}
\setminus\Gamma_{\mathcal{PZ}}\) the matrices \(F_{\alpha}(\boldsymbol{t})\) and
\(G_{\alpha}(\boldsymbol{t})\) are, respectively, the left zero and right pole
semiresidual matrices of \(R(\cdot,\boldsymbol{t})\):
\begin{equation}\label{kokopz}
F_{\alpha}(\boldsymbol{t})=F_{\mathcal{Z},\alpha}(\boldsymbol{t}),
\quad G_{\alpha}(\boldsymbol{t})=G_{\mathcal{P},\alpha}(\boldsymbol{t}),\quad
\forall\boldsymbol{t}\in\mathcal{D}_{\alpha}\setminus\Gamma_{\mathcal{PZ}}.
\end{equation}
It   is locally given by
\begin{equation}\label{jointrep2t}
R(z,\boldsymbol{t}) =I+F_{\alpha}(\boldsymbol{t})S_\alpha^{-1}(\boldsymbol{t})
(zI-A_{\mathcal{P}}(\boldsymbol{t}))^{-1}G_{\alpha}(\boldsymbol{t}),
\quad\boldsymbol{t}\in\mathcal{D}_{\alpha}\setminus\Gamma_{\mathcal{PZ}},\end{equation}
where \(\{S_\alpha(\boldsymbol{t})\}_{\alpha\in\mathfrak{A}}\) is the collection of
\(\mathcal{PZ}\)-coupling functions related to the pair of collections
\(\{F_{\alpha}(\boldsymbol{t})\}_{\alpha\in\mathfrak{A}}\),
\(\{G_{\alpha}(\boldsymbol{t})\}_{\alpha\in\mathfrak{A}}\). Furthermore, the
logarithmic derivative of \(R(\cdot,\boldsymbol{t})\) admits the local
representation
\begin{multline}
\dfrac{\partial R(z,\boldsymbol{t})}{\partial z}R^{-1}(z,\boldsymbol{t})=\\ =
-F_{\alpha}(\boldsymbol{t})S_\alpha^{-1}(\boldsymbol{t})
 (zI-A_{\mathcal{P}}(\boldsymbol{t}))^{-1}S_\alpha(\boldsymbol{t})
(zI-A_{\mathcal{Z}}(\boldsymbol{t}))^{-1}S_\alpha^{-1}(\boldsymbol{t})
G_{\alpha}(\boldsymbol{t}),\\\boldsymbol{t}\in\mathcal{D}_{\alpha}\setminus\Gamma_{\mathcal{PZ}}.
\label{replogder2t}
\end{multline}
\end{theorem}
\begin{proof} The proof is analogous to that of Theorem \ref{main1}.
\end{proof}

\section{Isosemiresidual families of generic rational matrix
 functions}
 \label{IFRMFGP}
In the present section we shall consider  an important special case
of holomorphic families of generic rational functions, parameterized
by the pole and zero loci \(\boldsymbol{t}\). Namely,
we are interested in the
case when (either left pole and right zero or right pole and left zero)
semiresidual functions of  \(\boldsymbol{t}\),
determining the family as explained in
Section \ref{HFGRF}, are constant.
\begin{definition}
\label{IsoSemiRes}
Let \(\boldsymbol{\mathcal{D}}\) be a domain in \(\mathbb{C}^{2n}_*\) and let
\(\{R(z,\boldsymbol{t})\}_{\boldsymbol{t}\in\boldsymbol{\mathcal{D}}}\) be a
normalized
holomorphic family of \(\mathbb{C}^{m\times m}\)-valued generic rational
functions, parameterized by the pole and zero loci.
\begin{enumerate}
\item
The family
\(\{R(z,\boldsymbol{t})\}_{\boldsymbol{t}\in\boldsymbol{\mathcal{D}}}\)
is said to be \(\mathcal{ZP}\)-{\em isosemiresidual}\footnote{Iso- (from \(\stackrel{.\!.}{\iota}\)\(\sigma o\varsigma\)
 \,-\,equal\,-\, in Old Greek)
 is a combining form.} if
there exists a pair of  matrices \(F\in\mathbb{C}^{m\times n}\) and
\(G\in\mathbb{C}^{n\times m}\) such that for every
\(\boldsymbol{t}\in\boldsymbol{\mathcal{D}}\) the matrices
\(F\) and \(G\) are, respectively,  the left pole and right zero
semiresidual matrices of the generic rational function
\(R(\cdot,\boldsymbol{t})\):
\[F=F_{\mathcal{P}}(\boldsymbol{t}),\quad G=G_{\mathcal{Z}}(\boldsymbol{t}),
\quad \forall \boldsymbol{t}\in\boldsymbol{\mathcal{D}}.\]

\item
The family
\(\{R(z,\boldsymbol{t})\}_{\boldsymbol{t}\in\boldsymbol{\mathcal{D}}}\)
is said to be \(\mathcal{PZ}\)-{\em isosemiresidual} if
there exists a pair of  matrices \(F\in\mathbb{C}^{m\times n}\) and
\(G\in\mathbb{C}^{n\times m}\) such that for every
\(\boldsymbol{t}\in\boldsymbol{\mathcal{D}}\) the matrices
\(F\) and \(G\) are, respectively,  the left zero and right pole
semiresidual matrices of the generic rational function
\(R(\cdot,\boldsymbol{t})\):
\[F=F_{\mathcal{Z}}(\boldsymbol{t}),\quad G=G_{\mathcal{P}}(\boldsymbol{t}),
\quad \forall \boldsymbol{t}\in\boldsymbol{\mathcal{D}}.\]
\end{enumerate}
\end{definition}

Let us assume that a pair of matrices \(F\in\mathbb{C}^{m\times n}\) and
\(G\in\mathbb{C}^{n\times m}\) is given. How to construct a
(\(\mathcal{PZ}\)- or \(\mathcal{ZP}\)-) isosemiresidual
normalized
holomorphic family of \(\mathbb{C}^{m\times m}\)-valued generic rational
functions, parameterized by the pole and zero loci, for which
the
 constant  functions
 \begin{equation}\label{consem}
F(\boldsymbol{t})\equiv F,\quad G(\boldsymbol{t})\equiv G
\end{equation}
would be the appropriate semiresidual functions? This  is a special case of the
problem  considered in Section  \ref{HFPD} (see Theorem \ref{main1}). Note, however,
that in this case the prescribed semiresidual functions \eqref{consem} are
holomorphic in the  domain \(\mathbb{C}^{2n}_*\). Therefore, we may consider its
open covering consisting of the single set -- the domain itself.

 Following the approach described in Section  \ref{HFPD},
we  consider
the solutions  \(S_{\mathcal{PZ}}(\boldsymbol{t})\),
\(S_{\mathcal{ZP}}(\boldsymbol{t})\) of the Lyapunov equations
\begin{subequations}\label{SLTP}
\begin{align}
\label{SL1TPrC}
A_{\mathcal{P}}(\boldsymbol{t})S_{\mathcal{PZ}}(\boldsymbol{t})-
S_{\mathcal{PZ}}(\boldsymbol{t})
A_{\mathcal{Z}}(\boldsymbol{t})&=
G F,\\
\label{SL2TPrC}
A_{\mathcal{Z}}(\boldsymbol{t})S_{\mathcal{ZP}}(\boldsymbol{t})-
S_{\mathcal{ZP}}(\boldsymbol{t})
A_{\mathcal{P}}(\boldsymbol{t})&=
G F,
\end{align}
\end{subequations}
where
\begin{equation}
A_{\mathcal{P}}(\boldsymbol{t})=\diag(t_1, \dots ,t_n), \quad
A_{\mathcal{Z}}(\boldsymbol{t})=\diag(t_{n+1}, \dots, t_{2n}). \label{polematrt1c}
\end{equation}
Then the functions \(S_{\mathcal{PZ}}(\boldsymbol{t})\),
\(S_{\mathcal{ZP}}(\boldsymbol{t})\) are given explicitly
by
\begin{subequations}
\label{Sol1C}%
\begin{align}\label{Sol1CPZ}
S_{\mathcal{PZ}}(\boldsymbol{t})&=
\begin{pmatrix}
\dfrac{g_{i} f_{j}}{t_{i}-t_{n+j}}
\end{pmatrix}_{1\leq i,j\leq n },\\
\label{Sol1CZP}
S_{\mathcal{ZP}}(\boldsymbol{t})&=
\begin{pmatrix}%
\dfrac{g_{i} f_{j}}{t_{n+i}-t_{j}}
\end{pmatrix}_{1\leq i,j\leq n },
\end{align}
\end{subequations}
where  \(g_i\) and \(f_{j}\) denote, respectively,
 the \(i\)-th row of \(G\)
and the \(j\)-th column of \(F\).
In particular, the functions \(S_{\mathcal{PZ}}(\boldsymbol{t})\),
\(S_{\mathcal{ZP}}(\boldsymbol{t})\) are  rational with respect to \(\boldsymbol{t}\)
and holomorphic in \(
\mathbb{C}^{2n}_{\ast}\).
The next step is to verify that the constant functions \eqref{consem}
are \(\mathcal{PZ}\)- or \(\mathcal{ZP}\)-admissible (that is, suitable for
the construction of a holomorphic family of generic rational functions -- see Theorem
\ref{main1}). This means to check that
\(
\det S_{\mathcal{PZ}}(\boldsymbol{t})\not\equiv 0\) or \(
\det S_{\mathcal{ZP}}(\boldsymbol{t})\not\equiv 0\). Note
that, since the functions \(S_{\mathcal{PZ}}(\boldsymbol{t})\),
\(S_{\mathcal{ZP}}(\boldsymbol{t})\) are identical up to the permutation
of variables \(t_k\leftrightarrow t_{n+k}\), \(1\leq k\leq n\),
these conditions are equivalent.
\begin{definition}
A pair of matrices \(F\in\mathbb{C}^{m\times n}\) and
\(G\in\mathbb{C}^{n\times m}\) is said to be {\em admissible}
if the \(\mathbb{C}^{n\times n}\)-valued rational
 functions \(S_{\mathcal{PZ}}(\boldsymbol{t})\),
\(S_{\mathcal{ZP}}(\boldsymbol{t})\) given by \eqref{Sol1C}
satisfy the (equivalent) conditions
\[
\det S_{\mathcal{PZ}}(\boldsymbol{t})\not\equiv 0,\quad
\det S_{\mathcal{ZP}}(\boldsymbol{t})\not\equiv 0.\]
\end{definition}

It turns out that the admissibility of a given  pair of matrices
\(F\in\mathbb{C}^{m\times n}\), \(G\in\mathbb{C}^{n\times m}\)
can be checked
by means of  a simple criterion, described below.

 Recall that for a matrix
\(M=\big(m_{i,j}\big)_{i,j=1}^{ n}\),
\begin{equation}
\label{DetExp}%
 \det M=\sum_{\sigma}(-1)^\sigma
m_{1,\sigma(1)}\cdots m_{n,\sigma(n)},
\end{equation}
where \(\sigma\) runs over all \(n!\) permutations of the set
\(1,\dots,n\), and \((-1)^\sigma\) is equal to either \(1\)
or \(-1\) depending on the parity of the permutation \(\sigma\).

\begin{definition}%
\label{DefZSingMatr}%
A matrix \(M\in\mathbb{C}^{n\times n}\) is said to be {\em Frobenius-singular}
if for some \(\ell,\) \(1\leq \ell\leq n,\)  there exist
indices \(1\leq\alpha_1<\dots<\alpha_\ell\leq n\);
\(1\leq\beta_1<\dots<\beta_{n-\ell+1}\leq n\),
 such that \(m_{\alpha_i,\beta_j}=0\) for all
\(1\leq i\leq \ell,1\leq j\leq n-\ell+1\).
\end{definition}%

\begin{theorem}%
\label{FrobTh}%
A matrix \(M\in\mathbb{C}^{n\times n}\) is
Frobenius-singular if and only if all \(n!\) summands
\((-1)^\sigma
m_{1,\sigma(1)}\cdots m_{n,\sigma(n)}\)
 of the
sum \eqref{DetExp} representing the determinant \(\det M\) are equal to zero.
\end{theorem}%

Theorem\ref{FrobTh}  is due to G.Frobenius, \cite{Fro1}. The proof
of this theorem can be also found in \cite{Ber}, Chapter 10, Theorem
9. The book \cite{LoPl} contains some historical remarks concerning
this theorem. See the Preface of \cite{LoPl}, especially pp.
xiii\,-\,xvii of the English edition (to which pp.\ 14\,-\,18 of the
Russian translation correspond).


\begin{proposition}%
\label{AdmCrit}%
A pair of matrices \(F\in\mathbb{C}^{m\times n}\),
\(G\in\mathbb{C}^{n\times m}\) is admissible if and only if their
product  \(GF\) is {\em not} a Frobenius-singular matrix.
\end{proposition}
\begin{proof}
Assume first that
the matrix \(GF\) is  Frobenius-singular.
 Then, according to \eqref{Sol1CPZ},
\begin{equation}
\label{DetCoupMatr}
\det S_{\mathcal{PZ}}(\boldsymbol{t})=
\sum_{\sigma}(-1)^\sigma\frac{
m_{1,\sigma(1)}\cdots m_{n,\sigma(n)}}{(t_1-t_{1+\sigma(1)})\cdots
((t_n-t_{n+\sigma(n)})},
\end{equation}
where \(\sigma\) runs over all \(n!\) permutations of the set
\(1,\dots,n\) and
\(m_{i,j}=g_{i}f_{n+j}\).
If  the matrix \(GF\) is Frobenius-singular  then,
according to Theorem
\ref{FrobTh}, all the numerators \(m_{1,\sigma(1)}\cdots m_{n,\sigma(n)}\)
 of the summands in \eqref{DetCoupMatr} are equal to
zero, and hence \(\det S_{\mathcal{PZ}}(\boldsymbol{t})\equiv 0\).

 Conversely,
if the matrix \(GF\) is not Frobenius-singular then,
 according to the
same Theorem \ref{FrobTh}, there exists a permutation \(\sigma_0\)
such that
\begin{equation*}%
m_{1,\sigma_0(1)}\cdots m_{n,\sigma_0(n)}\not=0.
\end{equation*}%
Let us choose and fix \(n\) pairwise different numbers
\(t_1^0,\,\dots\,\,t_n^0\) and set
\(t_{n+\sigma_0(1)}^0=t_1^0-\varepsilon,
\dots,t_{n+\sigma_0(n)}=t_n^0-\varepsilon\), where
\(\varepsilon\not=0\). Then as \(\varepsilon\to
0\) the summand
\begin{equation*}
(-1)^{\sigma_0}\frac{m_{1,\sigma_0(1)}\cdots
m_{n,\sigma_0(n)}}{(t_1-t_{n+\sigma_0(1)}(\varepsilon))
\cdots(t_n-t_{n+\sigma_0(n)}(\varepsilon))}
=(-1)^{\sigma_0}
(m_{1,{\sigma_0(1)}}\cdots
m_{n,{\sigma_0(n)}})\varepsilon^{-n}
\end{equation*}
is the leading term of the sum on the right-hand side of
\eqref{DetCoupMatr}: all other summands grow at most
as \(O(\epsilon^{-(n-1)})\).
\end{proof}
\begin{definition}
Let  a pair of matrices \(F\in\mathbb{C}^{m\times n}\),
\(G\in\mathbb{C}^{n\times m}\) be such that  their product
\(GF\) is {\em not}
 a Frobenius-singular matrix and let the
\(\mathbb{C}^{n\times n}\)-valued rational
 functions \(S_{\mathcal{PZ}}(\boldsymbol{t})\),
\(S_{\mathcal{ZP}}(\boldsymbol{t})\) be given by \eqref{Sol1C}.
\begin{enumerate}
\item
The set \begin{equation}
\label{SingSetCPZ}
\Gamma_{\mathcal{PZ}}=\{\boldsymbol{t}\in\mathbb{C}^{2n}_{\,\ast}:
\det S_{\mathcal{PZ}}(\boldsymbol{t})=0\}
\end{equation}
is said to be the \(\mathcal{PZ}\)-{\em singular set} related
to the pair \(F,G\).

\item
The set \begin{equation}
\label{SingSetCZP}
\Gamma_{\mathcal{ZP}}=\{\boldsymbol{t}\in\mathbb{C}^{2n}_{\,\ast}:
\det S_{\mathcal{ZP}}(\boldsymbol{t})=0\}
\end{equation}
 is said to be the \(\mathcal{ZP}\)-{\em singular set} related
to the pair \(F,G\).
\end{enumerate}
\end{definition}
\begin{remark}
Note that, since \(\det S_{\mathcal{PZ}}(\boldsymbol{t})\),
\(\det S_{\mathcal{ZP}}(\boldsymbol{t})\)
are  polynomials in
\((t_i-t_j)^{-1}\),
the singular sets \(\Gamma_{\mathcal{PZ}}\), \(\Gamma_{\mathcal{ZP}}\)
related
to the pair \(F,G\)  are  complex algebraic varieties of
codimension one in \(\mathbb{C}^{2n}_{\ast}\).
\end{remark}

Combining Proposition \ref{AdmCrit} and Theorems \ref{main1}, \ref{main2}, we obtain
\begin{theorem}
\label{Firmar}
Let  matrices \(F\in\mathbb{C}^{m\times n}\) and
\(G\in\mathbb{C}^{n\times m}\) be such that  their product
\(GF\) is {\em not}
 a Frobenius-singular matrix, and let
\(\Gamma_{\mathcal{PZ}}\), \(\Gamma_{\mathcal{ZP}}\) be the related singular sets.
Then the following statements hold true.
\begin{enumerate}

\item
There exists a unique \(\mathcal{PZ}\)-isosemiresidual family
\(\{R(z,\boldsymbol{t})\}_{\boldsymbol{t}
\in\mathbb{C}^{2n}_{\ast}\setminus \Gamma_{\mathcal{PZ}}}\)
of \(\mathbb{C}^{m\times m}\)-valued generic rational functions
 such that
for every
 \(\boldsymbol{t}\in\mathbb{C}^{2n}_{\ast}\setminus\Gamma_{\mathcal{PZ}}\)
the matrices \(F\) and \(G\)
are, respectively, the left zero and right pole semiresidual matrices
of \(R(\cdot,\boldsymbol{t})\):
\begin{equation}\label{kokopzt}
F=F_{\mathcal{Z}}(\boldsymbol{t}),
\quad G=G_{\mathcal{P}}(\boldsymbol{t}),\quad
\forall\boldsymbol{t}\in\mathbb{C}^{2n}_{\ast}\setminus\Gamma_{\mathcal{PZ}}.
\end{equation}
It is given by
\begin{equation}\label{joic2}
R(z,\boldsymbol{t}) =I+FS_{\mathcal{PZ}}^{-1}(\boldsymbol{t})
(zI-A_{\mathcal{P}}(\boldsymbol{t}))^{-1} G,\quad
\boldsymbol{t}\in\mathbb{C}^{2n}_{\ast}\setminus\Gamma_{\mathcal{PZ}}
\end{equation} where the  function
\(S_{\mathcal{PZ}}(\boldsymbol{t})\) satisfying the  equation \eqref{SL1TPrC} is
given by \eqref{Sol1CPZ}. Furthermore, the logarithmic derivative of
\(R(z,\boldsymbol{t})\) with respect to~\(z\) admits the representation
\begin{multline}
\dfrac{\partial R(z,\boldsymbol{t})}{\partial z}R^{-1}(z,\boldsymbol{t})=\\ =
-FS_{\mathcal{PZ}}^{-1}(\boldsymbol{t})
 (zI-A_{\mathcal{P}}(\boldsymbol{t}))^{-1}S_{\mathcal{PZ}}(\boldsymbol{t})
(zI-A_{\mathcal{Z}}(\boldsymbol{t}))^{-1}S_{\mathcal{PZ}}^{-1}(\boldsymbol{t})G,\\
\boldsymbol{t}\in\mathbb{C}^{2n}_{\ast}\setminus\Gamma_{\mathcal{PZ}}.
\label{replogc2t}
\end{multline}

\item  There exists a unique \(\mathcal{ZP}\)-isosemiresidual family
\(\{R(z,\boldsymbol{t})\}_{\boldsymbol{t}
\in\mathbb{C}^{2n}_{\ast}\setminus \Gamma_{\mathcal{ZP}}}\)
of \(\mathbb{C}^{m\times m}\)-valued generic rational functions
 such that
for every
 \(\boldsymbol{t}\in\mathbb{C}^{2n}_{\ast}\setminus\Gamma_{\mathcal{ZP}}\)
the matrices \(F\) and \(G\)
are, respectively, the left pole and right zero semiresidual matrices
of \(R(\cdot,\boldsymbol{t})\):
\begin{equation}\label{kokozpt}
F=F_{\mathcal{P}}(\boldsymbol{t}),
\quad G=G_{\mathcal{Z}}(\boldsymbol{t}),\quad
\forall\boldsymbol{t}\in\mathbb{C}^{2n}_{\ast}\setminus\Gamma_{\mathcal{ZP}}.
\end{equation}
It is given by
\begin{equation}\label{joic1}
R(z,\boldsymbol{t}) =I-F
(zI-A_{\mathcal{P}}(\boldsymbol{t}))^{-1}S_{\mathcal{ZP}}^{-1}(\boldsymbol{t}) G,
\quad \boldsymbol{t}\in\mathbb{C}^{2n}_{\ast}\setminus\Gamma_{\mathcal{ZP}},
\end{equation} where the function
\(S_{\mathcal{ZP}}(\boldsymbol{t})\) satisfying the  equation \eqref{SL2TPrC} is
given by \eqref{Sol1CZP}. Furthermore, the logarithmic derivative of
\(R(z,\boldsymbol{t})\) with respect to~\(z\) admits the representation
\begin{multline}
\dfrac{\partial R(z,\boldsymbol{t})}{\partial z}R^{-1}(z,\boldsymbol{t})=\\ = F
 (zI-A_{\mathcal{P}}(\boldsymbol{t}))^{-1}S_{\mathcal{ZP}}^{-1}(\boldsymbol{t})
(zI-A_{\mathcal{Z}}(\boldsymbol{t}))^{-1}G,\\
\boldsymbol{t}\in\mathbb{C}^{2n}_{\ast}\setminus\Gamma_{\mathcal{ZP}}.
 \label{replogc1t}
\end{multline}
\end{enumerate}
\end{theorem}

\section{Isoprincipal families of generic rational matrix functions}
\label{IFGRMF2}

Our interest in holomorphic families of generic rational functions
is motivated by our intent to construct rational solutions
 of the Schlesinger system (see Section \ref{rsotss} below).
Indeed, given a holomorphic family
\(\{R(z,\boldsymbol{t}):\boldsymbol{t}\in\boldsymbol{\mathcal{D}}\}\) of
 generic rational  functions, we can consider
the linear differential system
\begin{equation}
\label{holde} \dfrac{\partial R(z,\boldsymbol{t})}{\partial
z}=Q_R(z,\boldsymbol{t})R(z,\boldsymbol{t}),
\end{equation}
where \(Q_R(z,\boldsymbol{t})\) is the logarithmic derivative of
\(R(z,\boldsymbol{t})\) with respect to \(z\). According to Lemma \ref{HRM} (and in
view of \eqref{parlogder}, \eqref{AELDT}), the system \eqref{holde} can be rewritten
as
\begin{equation}
\label{holdef} \dfrac{\partial R(z,\boldsymbol{t})}{\partial
z}=\left(\sum_{k=1}^{2n}\frac{Q_k(\boldsymbol{t})}{z-t_k}\right)R(z,\boldsymbol{t}),
\end{equation}
where the functions \(Q_k(\boldsymbol{t})\) are holomorphic in
\(\boldsymbol{\mathcal{D}}\). The system \eqref{holdef} can be
viewed as a holomorphic family (=deformation) of Fuchsian systems
parameterized by the singularities' loci. It was proved in
\cite{KaVo} (Theorem 8.2) that in the case when the deformation
\eqref{holdef} is {\em isoprincipal} the functions
\(Q_k(\boldsymbol{t})\)
 satisfy the Schlesinger system.
\begin{definition}%
\label{DefIsoPrinc}%
Let \(\boldsymbol{\mathcal{D}}\) be a domain in
\(\mathbb{C}^{2n}_*\) and let
\(\{R(z,\boldsymbol{t})\}_{\boldsymbol{t}\in\boldsymbol{\mathcal{D}}}\)
be a normalized holomorphic family of \(\mathbb{C}^{m\times
m}\)-valued generic rational functions, parameterized by the pole
and zero loci. Assume that for \(1\leq k\leq 2n\) there exist
\(\mathbb{C}^{m\times m}\)-valued functions \(E_k(\cdot)\),
holomorphic and invertible in \(\mathbb{C}_*\), such that for every
\(\boldsymbol{t}\in\boldsymbol{\mathcal{D}}\) the function \(E_k\)
is the principal factor\footnote{See Definition \ref{DefPrinc}.} of
the function \(R(\cdot,\boldsymbol{t})\) at \(t_k\): there exists a
\(\mathbb{C}^{m\times m}\)-valued function
\(H_k(\cdot,\boldsymbol{t})\),  holomorphic and invertible  in a
neighborhood of \(t_k\), such that
\begin{equation}\label{muldecot}
R(z,\boldsymbol{t})=H_k(z,\boldsymbol{t})E_k(z-t_k).
\end{equation}
 Then the  family
\(\{R(z,\boldsymbol{t})\}_{\boldsymbol{t}\in\boldsymbol{\mathcal{D}}}\) is said to
be {\em isoprincipal}.
\end{definition}%
\begin{theorem}
\label{CoincIso}%
Let \(\boldsymbol{\mathcal{D}}\) be a domain in
\(\mathbb{C}^{2n}_*\) and let
\(\{R(z,\boldsymbol{t})\}_{\boldsymbol{t}\in\boldsymbol{\mathcal{D}}}\)
be a normalized holomorphic family of \(\mathbb{C}^{m\times
m}\)-valued generic rational functions, parameterized by the pole
and zero loci. The family
\(\{R(z,\boldsymbol{t}):\boldsymbol{t}\in\boldsymbol{\mathcal{D}}\}\)
is isoprincipal if and only if it is
\(\mathcal{PZ}\)-isosemiresidual.
\end{theorem}
\begin{proof}
First, assume that
 the
family
\(\{R(z,\boldsymbol{t})\}_{\boldsymbol{t}\in\boldsymbol{\mathcal{D}}}\)
is \(\mathcal{PZ}\)-isosemiresidual. Then, according to Definition
\ref{IsoSemiRes},
 there exist
\(F\in\mathbb{C}^{m\times n}\) and \(G\in\mathbb{C}^{n\times m}\) such that for
every \(\boldsymbol{t}\in\boldsymbol{\mathcal{D}}\) the matrices \(F\) and \(G\)
are, respectively, the  left zero and right pole semiresidual matrices  of the
generic rational function \(R(\cdot,\boldsymbol{t})\). Then, by Lemma
\ref{RegPrinFact} (see also Remark \ref{semi2}), for \(k=1,\dots,2n\) and
{independently} of \(\boldsymbol{t}\) the principal factor \(E_k(\zeta)\) of
\(R(\cdot,\boldsymbol{t})\) at \(t_k\) can be chosen in the form
\begin{equation}\label{prinfac}
E_k(\zeta)=\left\{\begin{array}{l@{\quad\text{if}\quad}l}
I+L_k-\zeta^{-1}L_k,&1\leq k\leq n,\\[1ex]
I-L_k+\zeta L_k,&1+n\leq k\leq 2n,
\end{array}\right.
\end{equation}
where
\begin{equation}\label{expprin}
L_k=\left\{\begin{array}{l@{\quad\text{if}\quad}l}
-g_k^*(g_kg_k^*)^{-1}g_k,& 1\leq k\leq n,\\[1ex]
\phantom{-}f_k(f_k^*f_k)^{-1}f_k^*,& 1+n\leq k\leq 2n,
\end{array}\right.
\end{equation}
\(g_k\) is the \(k\)-th row of \(G\) and \(f_{k}\) is the
\((k-n)\)-th column of \(F\). Hence, by Definition \ref{DefIsoPrinc}
the family
\(\{R(z,\boldsymbol{t}):\boldsymbol{t}\in\boldsymbol{\mathcal{D}}\}\)
is isoprincipal.

Conversely, assume that the family
\(\{R(z,\boldsymbol{t})\}_{\boldsymbol{t}\in\boldsymbol{\mathcal{D}}}\) is
isoprincipal. Then, according to Definition \ref{DefIsoPrinc}, for \(1\leq k\leq
2n\) there exist functions
 \(\hat{E_k}(\cdot)\), holomorphic and invertible in \(\mathbb{C}_*\), such that for every
\(\boldsymbol{t}\in\boldsymbol{\mathcal{D}}\) the function \(\hat{E_k}\) is the
principal factor of the function \(R(\cdot,\boldsymbol{t})\) at \(t_k\). Let
\(\boldsymbol{t}^0,\boldsymbol{t}\in\boldsymbol{\mathcal{D}}\) be fixed and let
\(F\in\mathbb{C}^{m\times n}\) and \(G\in\mathbb{C}^{n\times m}\) be, respectively,
the  left zero and right pole semiresidual matrices  of the generic rational
function \(R(\cdot,\boldsymbol{t}^0)\). Let us denote by \(g_k\)  the \(k\)-th row
of \(G\) and by \(f_{k}\)  the \((k-n)\)-th column of \(F\). Then, in view of Remark
\ref{gauprr},
 the
function \(\hat{E_k}\) is of the form \[ \hat{E_k}(\zeta)=M_k(\zeta)E_k(\zeta),\]
where \(M_k(\cdot)\) is a \(\mathbb{C}^{m\times m}\)-valued function, holomorphic
and invertible in \(\mathbb{C}\) and \(E_k(\cdot)\) is
 given by \eqref{prinfac}, \eqref{expprin}. Hence for \(z\) in a neighborhood of  \(t_k\)
\(R(z,\boldsymbol{t})\) admits the representation
\[R(z,\boldsymbol{t})=H_k(z,\boldsymbol{t})M_k(z-t_k)E_k(z-t_k),\]
where \(H_k(\cdot,\boldsymbol{t})\) is a \(\mathbb{C}^{m\times m}\)-valued
function,
holomorphic and invertible at \(t_k\). Then, for \(k=1,\ldots,n\), the residue
\(R_k(\boldsymbol{t})\) of \(R(z,\boldsymbol{t})\) at \(t_k\)  is given by
\[R_k(\boldsymbol{t})=\left(H_k(t_k,\boldsymbol{t})M_k(0)g_k^*(g_kg_k^*)^{-1}\right)g_k.\]
Therefore, \(G\) is the right pole semiresidual matrix
of the  function
\(R(z,\boldsymbol{t})\), as well. Analogously,
for \(k=n+1,\ldots,2n\)
\[E^{-1}(\zeta)=I-L_k+\zeta^{-1}L_k,\]
where
\[L_k=f_k(f_k^*f_k)^{-1}f_k^*,\]
hence the residue \(R_k(\boldsymbol{t})\)
of \(R^{-1}(z,\boldsymbol{t})\) at \(t_k\)
is given by
\[R_k(\boldsymbol{t})=
f_k\left((f_k^*f_k)^{-1}f_k^*M_k^{-1}(0)H_k^{-1}(t_k,\boldsymbol{t})\right).\]
 Therefore, \(F\) is the left zero semiresidual matrix
of the  function
\(R(z,\boldsymbol{t})\), as well.
This completes the proof.
\end{proof}
Theorem \ref{CoincIso} reduces  the construction of an {\em
isoprincipal} family  to the
 construction of an {\em isosemiresidual} family.
The latter problem has already been considered in Section \ref{IFRMFGP}. According
to Theorems \ref{Firmar} and \ref{CoincIso}, from any  pair of matrices
\(F\in\mathbb{C}^{m\times n}\) and \(G\in\mathbb{C}^{n\times m}\), such that the
product \(GF\) is { not}
 a Frobenius-singular matrix, we can construct an isoprincipal family
of generic rational functions
\(\{R(z,\boldsymbol{t})\}_{\boldsymbol{t}\in\mathbb{C}^{2n}_*\setminus\Gamma_{\mathcal{PZ}}}\),
where \(\Gamma_{\mathcal{PZ}}\) denotes the \(\mathcal{PZ}\)-singular set related to
the pair of matrices \(F\), \(G\). This family is given by
\begin{equation}\label{joic2o2}
R(z,\boldsymbol{t}) =I+FS_{\mathcal{PZ}}^{-1}(\boldsymbol{t})
(zI-A_{\mathcal{P}}(\boldsymbol{t}))^{-1} G, \end{equation} where the function
\(S_{\mathcal{PZ}}(\boldsymbol{t})\), satisfying \eqref{SL1TPrC}, is given by
\eqref{Sol1CPZ}. The logarithmic derivative of \(R(z,\boldsymbol{t})\) with respect
to \(z\) is given by
 \begin{multline}\label{almost1}
\dfrac{\partial R(z,\boldsymbol{t})}{\partial z}R^{-1}(z,\boldsymbol{t})=\\ =
-FS_{\mathcal{PZ}}^{-1}(\boldsymbol{t})
 (zI-A_{\mathcal{P}}(\boldsymbol{t}))^{-1}S_{\mathcal{PZ}}(\boldsymbol{t})
(zI-A_{\mathcal{Z}}(\boldsymbol{t}))^{-1}S_{\mathcal{PZ}}^{-1}(\boldsymbol{t})G,
\end{multline}
 and we obtain the following expressions for its residues \(Q_k(\boldsymbol{t})\)
\begin{subequations}\label{almost2}
\begin{multline}
Q_k(\boldsymbol{t})= -FS_{\mathcal{PZ}}^{-1}(\boldsymbol{t})
 I_{[k]}S_{\mathcal{PZ}}(\boldsymbol{t})
(t_k I-A_{\mathcal{Z}}(\boldsymbol{t}))^{-1}
S_{\mathcal{PZ}}^{-1}(\boldsymbol{t})G,\\ 1\leq k\leq n,\end{multline}
\begin{multline}
Q_k(\boldsymbol{t})=-FS_{\mathcal{PZ}}^{-1}(\boldsymbol{t})
 (t_kI-A_{\mathcal{P}}(\boldsymbol{t}))^{-1}S_{\mathcal{PZ}}(\boldsymbol{t})
I_{[k-n]}S_{\mathcal{PZ}}^{-1}(\boldsymbol{t})G,\\ n+1\leq k\leq 2n.
\end{multline}
\end{subequations}
Here we use the notation
\[I_{[k]}\stackrel{\textup{\tiny def}}{=}
\diag(\delta_{1,k},\dots,\delta_{n,k}),\]
where \(\delta_{i,j}\) is the Kronecker delta.

\begin{remark}
Note that, according to \eqref{Sol1CPZ},  the function
\(S_{\mathcal{PZ}}(\boldsymbol{t})\) is a rational function
of \(\boldsymbol{t}\). Hence also the functions
\(Q_k(\boldsymbol{t})\) are rational functions
of \(\boldsymbol{t}\).
\end{remark}
\section{Rational solutions of the Schlesinger system}
\label{rsotss}
 It can be checked
  that the rational functions \(Q_k(\boldsymbol{t})\) given by
\eqref{almost2} satisfy the Schlesinger system
\begin{equation}
\label{tscs}
\left\{\begin{array}{rl} \dfrac{\partial Q_k}{\partial
t_\ell}&=\dfrac{[Q_\ell,Q_k]}{t_\ell-t_k},
\quad k\not=\ell,\\[2ex]
\dfrac{\partial Q_k}{\partial t_k}&=\displaystyle\sum_{\ell\not=k}\dfrac{[Q_\ell,Q_k]}
{t_k-t_\ell}.
\end{array}
\right.
\end{equation}
 It is also not very difficult to check that
\begin{equation}%
\label{PotFuncNP}%
V(\boldsymbol{t})\stackrel{\textup{\tiny def}}{=}
-FS_{\mathcal{PZ}}^{-1}(\boldsymbol{t})G
\end{equation}%
is the potential function for this solution:
\begin{equation}%
\label{PotNP}%
Q_k(\boldsymbol{t})=\frac{\partial V(\boldsymbol{t})} {\partial t_k},\quad
k=1,\dots,2n.
\end{equation}%

Furthermore, one can show that the rational function \(\det
S_{\mathcal{PZ}}(\boldsymbol{t})\) admits the following integral representation
\begin{equation}
\det S_{\mathcal{PZ}}(\boldsymbol{t})= \det
S_{\mathcal{PZ}}(\boldsymbol{t}_0)\cdot
\exp\bigg\{\int\limits_{\gamma} \sum\limits_{\substack{1\leq i,j\leq
2n,\\j\not=i}}\dfrac{\trace\Big( \frac{\partial
V(\boldsymbol{t})}{\partial t_i}\cdot \frac{\partial
V(\boldsymbol{t})}{\partial t_j}\Big)}{t_i-t_j}\,\,dt_i\bigg\}\,.
\end{equation}
where \(\boldsymbol{t}_0\) and \(\boldsymbol{t}\) are two arbitrary points the
domain \(\mathbb{C}^{2n}_{\,\ast}\setminus\Gamma_{\mathcal{PZ}}\), and \(\gamma\) is
an arbitrary path which begins at \(\boldsymbol{t}_0\), ends at \(\boldsymbol{t}\)
and is contained in \(\mathbb{C}^{2n}_{\,\ast}\setminus\Gamma_{\mathcal{PZ}}\).

However, the explanation of these facts lies in the considerations of Sections 2 and
3 of the first part \textup{\cite{KaVo}} of this work. The matrix functions
\(Q_k(\boldsymbol{t})\) satisfy the Schlesinger system, and the function
\(-V(\boldsymbol{t})\) is a Laurent coefficient at \(z=\infty\) of the normalized
solution \eqref{joic2o2} of the Fuchsian system
\begin{equation}%
\label{FuchSy}%
\frac{d\,R(z,\boldsymbol{t})}{d\,z}= \left( \sum\limits_{1\leq k\leq
2n}\frac{Q_k(\boldsymbol{t})}{z-t_k} \right) R(z,\boldsymbol{t}),
\end{equation}%
\begin{equation}%
\label{ResAtInf}%
R(z,\boldsymbol{t})=I-\frac{V
(\boldsymbol{t})}{z}+o(|z|^{-1})\,\,\textup{ as}\ %
\ z\to\infty,
\end{equation}%
while the function
\begin{equation}
\tau(\boldsymbol{t})\stackrel{\textup{\tiny def}}{=}\det
S_{\mathcal{PZ}}(\boldsymbol{t}),
\end{equation}
is the tau-function related to the solution
\(Q_1(\boldsymbol{t}),\dots,Q_{2n}(\boldsymbol{t})\) of the Schlesinger system.

More detailed explanation of these and other related facts will be
given in the third part of this work.\\[2ex]

\renewcommand{\thesection}{\Alph{section}}
\setcounter{section}{1}%
\[ \boldsymbol{\mathcal{APPENDIX}}\]

\section{The global factorization of a holomorphic matrix function of rank
one}\label{app2}

Let \(M(\boldsymbol{t})=\|m_{p,q}(\boldsymbol{t})\|_{1\leq p,q \leq
m}\) be a \(\mathbb{C}^{m\times m}\)-valued  function of the
variable \(\boldsymbol{t}\in\boldsymbol{\mathcal{D}}\), where
\(\boldsymbol{\mathcal{D}}\)  is a domain in \(\mathbb{C}^N\). (We
can even assume that \(\boldsymbol{\mathcal{D}}\) is a Riemann
domain\footnote{See Definition 5.4.4 in \cite{Her}.} of dimension
\(N\) over \(\mathbb{C}^N\).)  In our considerations \(N=2n\) and
\(\boldsymbol{\mathcal{D}}\subseteq\mathbb{C}^{2n}_{\ast}.\)
 Let the matrix function \(M\)
be holomorphic in \(\boldsymbol{\mathcal{D}}\) and let
\begin{equation}\label{RONE}
\rank M(\boldsymbol{t})=1\quad\forall
\boldsymbol{t}\in\boldsymbol{\mathcal{D}}. \end{equation}
 We will
try to represent \(M\) in the form
\begin{equation}%
\label{DesirFact}%
M(\boldsymbol{t})=f(\boldsymbol{t})g(\boldsymbol{t}),
\end{equation}%
where\(f(\boldsymbol{t})\) and \(g(\boldsymbol{t})\) are,
respectively, a \(\mathbb{C}^{m\times 1}\)-valued function and a
\(\mathbb{C}^{1\times m}\)-valued function, both of them
holomorphic\footnote{In general, such a global factorization is
impossible even if the factors \(f(\boldsymbol{t})\) and
\(g(\boldsymbol{t})\) are only required to be  continuous rather
than holomorphic: one of the obstacles is of topological nature.} in
\(\boldsymbol{\mathcal{D}}\).

Let us recall that, according to Lemma \ref{RankOneMaLe},
 there exist a finite open covering \(\{\mathcal{U}_p\}_{p=1}^m\) %
of  \(\boldsymbol{\mathcal{D}}\), a collection
\(\{f_p(\boldsymbol{t})\}_{p=1}^m\) of \(\mathbb{C}^{m\times
1}\)-valued functions and a collection \(
\{g_p(\boldsymbol{t})\}_{p=1}^m\) of \(\mathbb{C}^{1\times
m}\)-valued functions satisfying the following conditions.
\begin{enumerate}
\item
For \(p=1,\dots,m\) the  functions \(f_p(\boldsymbol{t})\) and
\(g_p(\boldsymbol{t})\) are holomorphic in \(\mathcal{U}_p\).

\item For \(p=1,\dots,m\) the function
\(M(\boldsymbol{t})\) admits the factorization
\begin{equation}%
\label{LocFactB}%
M(\boldsymbol{t})=f_p(\boldsymbol{t})g_p(\boldsymbol{t}),\quad
\boldsymbol{t}\in\boldsymbol{\mathcal{U}}_p.
\end{equation}%
\item
Whenever
\(\mathcal{U}_{p^\prime}\cap\mathcal{U}_{p^{\prime\prime}}\not=\emptyset\),
there exists a (scalar) function
\(\varphi_{p^\prime,p^{\prime\prime}}(\boldsymbol{t})\), holomorphic
and invertible in
\(\mathcal{U}_{p^\prime}\cap\mathcal{U}_{p^{\prime\prime}}\), such
that
\begin{equation}\label{gauauB}
f_{p^{\prime\prime}}(\boldsymbol{t})=f_{p^{\prime}}(\boldsymbol{t})
\varphi_{p^\prime,p^{\prime\prime}}(\boldsymbol{t}),\
g_{p^{\prime\prime}}(\boldsymbol{t})=
\varphi^{-1}_{p^\prime,p^{\prime\prime}}(\boldsymbol{t})g_{p^{\prime}}(\boldsymbol{t})
\quad\forall\boldsymbol{t}\in\mathcal{U}_{p^\prime}\cap\mathcal{U}_{p^{\prime\prime}}.
\end{equation}
In particular,
\begin{subequations}%
\label{CocycPropAp}%
\begin{align}
\varphi_{p,p}(\boldsymbol{t})=1 &\quad \forall \boldsymbol{t}\in
\mathcal{U}_p, \\
\varphi_{p,p^\prime}(\boldsymbol{t})=\varphi_{p,p^{\prime\prime}}(\boldsymbol{t})
\varphi_{p^{\prime\prime},p^\prime}(\boldsymbol{t}) &\quad \forall
\boldsymbol{t}\in \mathcal{U}_p\cap
\mathcal{U}_{p^\prime}\cap\mathcal{U}_{p^{\prime\prime}}.
\end{align}
\end{subequations}%
\end{enumerate}

The equalities \eqref{LocFactB}, \(p=1,\,\dots\,,\,k\), are nothing
more than the factorizations of the form \eqref{DesirFact}, with
holomorphic factors \(f_p(\boldsymbol{t})\) and
\(g_p(\boldsymbol{t})\). However, the factorization \eqref{LocFactB}
is only {\em local}: for each \(p\) the equality \eqref{LocFactB}
holds in the open subset \(\mathcal{U}_p\) of the set
\(\boldsymbol{\mathcal{D}}\). For different \(p^\prime\) and
\(p^{\prime\prime}\), the factors \(f_{p^\prime}\,, g_{p^\prime}\)
and \(f_{p^{\prime\prime}}\,, g_{p^{\prime\prime}}\) may {\em not
agree in the intersections} \(\mathcal{U}_{p^\prime}\cap
\mathcal{U}_{p^{\prime\prime}}\). To glue  the factorizations
\eqref{LocFactB} for different \(p\) together, we seek scalar
functions \(\varphi_p(\boldsymbol{t})\) which are holomorphic in
\(\mathcal{U}_p\), do not  vanish there  and satisfy the condition
\begin{equation}%
\label{AgrCond1}%
f_{p^\prime}(\boldsymbol{t})\varphi_{p^\prime}(\boldsymbol{t}) =
f_{p^{\prime\prime}}(\boldsymbol{t})\varphi_{p^{\prime\prime}}(\boldsymbol{t})\quad
\forall\boldsymbol{t}\in \mathcal{U}_{p^\prime}\cap
\mathcal{U}_{p^{\prime\prime}}.
\end{equation}%
Then, in view of \eqref{LocFactB},
\begin{equation}%
\label{AgrCond2}%
\varphi_{p^\prime}^{-1}(\boldsymbol{t})
g_{p^\prime}(\boldsymbol{t})=\varphi_{p^{\prime\prime}}^{-1}(\boldsymbol{t})
g_{p^{\prime\prime}}(\boldsymbol{t})\quad \forall\boldsymbol{t}\in
\mathcal{U}_{p^\prime}\cap \mathcal{U}_{p^{\prime\prime}}.
\end{equation}%
Assuming that such functions  \(\varphi_{\,p}\,,\,%
1\leq p\leq k\), are found, we set
\begin{subequations}%
\label{DefFactors}%
\begin{align}
f(\boldsymbol{t})&\stackrel{\textup{\tiny def}}{=}
f_p(\boldsymbol{t})\varphi_p(\boldsymbol{t})\text{ if }
\boldsymbol{t}\in\mathcal{U}_p,\\
g(\boldsymbol{t})&\stackrel{\textup{\tiny
def}}{=}\varphi_{p}^{-1}(\boldsymbol{t}) g_p(\boldsymbol{t})\text{
if }
\boldsymbol{t}\in\mathcal{U}_p. %
\end{align}
\end{subequations}%
The relations \eqref{AgrCond1},  \eqref{AgrCond2} ensure that these
definitions are not contradictory. Thus the functions
\(f(\boldsymbol{t})\) and \(g(\boldsymbol{t})\) are defined for
every \(\boldsymbol{t}\in\boldsymbol{\mathcal{D}}\). Moreover, these
functions  are holomorphic in \(\boldsymbol{\mathcal{D}}\) and
provide the factorization \eqref{DesirFact}.

From \eqref{gauauB} it follows that the condition \eqref{AgrCond1}
is equivalent to the condition
\begin{equation}%
\label{CoboundAp}%
\varphi_{p^\prime}(\boldsymbol{t})
=\varphi_{p^\prime,p^{\prime\prime}}(\boldsymbol{t})
\varphi_{p^{\prime\prime}}(\boldsymbol{t})\quad
\forall\boldsymbol{t}\in \mathcal{U}_{p^\prime}\cap
\mathcal{U}_{p^{\prime\prime}},
\end{equation}%
where \(\varphi_{p^\prime,p^{\prime\prime}}(\boldsymbol{t})\) are
the functions appearing in \eqref{gauauB}. Thus, to ensure that the
 conditions \eqref{AgrCond1}, \eqref{AgrCond2} are in force, we have
to solve the so-called {\em second Cousin problem} (see \cite{Shab},
\cite{Oni}, \cite{Leit} and \cite{Her}):

\begin{Pb}\label{secop}
 Let \(\boldsymbol{\mathcal{D}}\) be a
complex manifold and let \(\{\mathcal{U}_{\alpha}\}_{\alpha}\) of
\(\boldsymbol{\mathcal{D}}\) be an open covering of
\(\boldsymbol{\mathcal{D}}\). For each \(\alpha,\,\beta\) such that
\(\mathcal{U}_{\alpha}\cap \mathcal{U}_{\beta}\not=\emptyset\) let a
\(\mathbb{C}\)-valued function \(\varphi_{\alpha,\beta}\),
holomorphic and non-vanishing in  \(\mathcal{U}_{\alpha}\cap
\mathcal{U}_{\beta}\), be given.

Find a collection of \(\mathbb{C}\)-valued functions
\(\{\varphi_{\alpha}\}_{\alpha}\) with the following properties:
\begin{enumerate}
\item
For every \(\alpha\) the function \(\varphi_{\alpha}\) is
holomorphic in \(\mathcal{U}_{\alpha}\) and does not vanish there.
\item
Whenever \(\mathcal{U}_{\alpha}\cap
\mathcal{U}_{\beta}\not=\emptyset\), the relation
\begin{equation}%
\label{CoboundCous}%
\varphi_{\alpha} =\varphi_{\alpha,\beta} \varphi_{\beta}
\end{equation}
holds in \(\mathcal{U}_{\alpha}\cap \mathcal{U}_{\beta}\).
\end{enumerate}
\end{Pb}

A necessary condition for the solvability of the second Cousin
problem in \(\boldsymbol{\mathcal{D}}\) with the given data
\(\{\mathcal{U}_{\alpha}\}_{\alpha},\,\{\varphi_{\alpha,\beta}\}_{\alpha,\beta}\)
is the so-called "{\em cocycle condition}":
\begin{equation}\label{CocycCondAp}
\begin{array}{rcll}%
 \varphi_{\alpha,\gamma}&=&\varphi_{\alpha,\beta}
\varphi_{\beta,\gamma}&\text{ in every non-empty triple intersection
}  \mathcal{U}_{\alpha}\cap \mathcal{U}_{\beta}\cap
\mathcal{U}_{\gamma},\\
 \varphi_{\alpha,\alpha}&= &1 &
\text{ in every }\mathcal{U}_{\alpha}.
\end{array}\end{equation}
In our case this condition is fulfilled -- see \eqref{CocycPropAp}.
However, the cocycle condition alone is not sufficient to guarantee
the existence of a solution to the second Cousin problem -- it
depends on \(\boldsymbol{\mathcal{D}}\) itself, as well.

\begin{proposition} \textup{(J.-P.Serre, \cite{Ser1}; see also
\cite{Shab}, section 16; \cite{Her},
 sections 5.5 and 7.4; \cite{Oni}, section 4.4.)}
If \(\boldsymbol{\mathcal{D}}\) is a Stein manifold\footnote{A
domain \(\boldsymbol{\mathcal{D}}\) in \(\mathbb{C}^N\) is a Stein
manifold if and only if \(\boldsymbol{\mathcal{D}}\) is
pseudoconvex, or, what is equivalent,  \(\boldsymbol{\mathcal{D}}\)
is a holomorphy domain.} which satisfies the
condition\,\footnote{\,\(H^2(\boldsymbol{\mathcal{D}},\,\mathbb{Z})\)
is the second cohomology group of \(\boldsymbol{\mathcal{D}}\) with
integer coefficients.}
\begin{equation}%
\label{HomTriv}%
 H^2(\boldsymbol{\mathcal{D}},\,\mathbb{Z})=0,
\end{equation}%
 then the second Cousin problem in \(\boldsymbol{\mathcal{D}}\)
with   arbitrary given data \(\{U_{\alpha},
\varphi_{\beta,\alpha}\}\) satisfying the cocycle condition
\eqref{CocycCondAp} is solvable.
\end{proposition}

As we have seen, the factorization problem \eqref{DesirFact} can be
reduced to solving the second Cousin with a certain data. Thus, the
following result holds:

\begin{theorem}%
\label{MatrFactRes}%
Let \(M(\boldsymbol{t})\) be a \(\mathbb{C}^{m\times m}\)-valued
function, holomorphic for
\(\boldsymbol{t}\in\boldsymbol{\mathcal{D}}\), where
\(\boldsymbol{\mathcal{D}}\) is a Riemann domain over
\(\mathbb{C}^N\). Assume that \(M\) satisfies the condition
\begin{equation*}
 \textup{rank}\,M(\boldsymbol{t})=1\quad \forall
 \boldsymbol{t}\in\boldsymbol{\mathcal{D}}.
\end{equation*}
 If \(\boldsymbol{\mathcal{D}}\) possesses the property:
{\em the second Cousin problem in \(\boldsymbol{\mathcal{D}}\) with
arbitrary given data
\(\{U_{\alpha}\}_\alpha,\,\{\varphi_{\alpha,\beta}\}_{\alpha,\beta}\)
satisfying the cocycle condition \eqref{CocycCondAp} is solvable,}
then the matrix function \(M(\boldsymbol{t})\) admits the
factorization of the form
\begin{equation*}%
 M(\boldsymbol{t})=f(\boldsymbol{t})\cdot
g(\boldsymbol{t}),
\end{equation*}%
where the factors \(f(\boldsymbol{t})\) and \(g(\boldsymbol{t})\)
are, respectively, a \(\mathbb{C}^{m\times 1}\)-valued function and
a \(\mathbb{C}^{1\times m}\)-valued function, holomorphic and
non-vanishing for \(\boldsymbol{t}\in\boldsymbol{\mathcal{D}}\). In
particular, such is the case if \(\boldsymbol{\mathcal{D}}\) is a
Stein manifold satisfying the condition \eqref{HomTriv}.
\end{theorem}


\end{document}